\documentclass[journal]{IEEEtran}
\usepackage{amsmath,amscd,amsbsy,amssymb,latexsym,url,bm,amsthm}
\allowdisplaybreaks
\usepackage{array}
\usepackage{algorithm}
\usepackage{algorithmic}
\usepackage{booktabs}
\usepackage{cite}
\usepackage{etoolbox} 
\usepackage[T1]{fontenc}
\usepackage[para]{footmisc} 
\usepackage{graphicx}
\usepackage{multirow,makecell}
\usepackage{subeqnarray}
\usepackage[caption=false,font=footnotesize,subrefformat=parens,labelformat=parens]{subfig}
\usepackage{tabularx}
\usepackage{threeparttable}
\usepackage{tikz}
\usepackage{xcolor}
\usepackage{balance}
\usepackage{hyperref}
\hypersetup{hidelinks}

\let\classAND\AND
\let\AND\relax
\usepackage{algorithmic}

\let\AND\classAND
\AtBeginEnvironment{algorithmic}{\let\AND\algoAND}

\newtheorem{assumption}{Assumption}

\newtheorem{theorem}{Theorem}
\newtheorem{lemma}{Lemma}
\newtheorem{remark}{Remark}

\newcommand{\bigO}[1]{\ensuremath{\mathop{}\mathopen{}\mathcal{O}\mathopen{}\left(#1\right)}}
\newcommand{\bigOcom}[1]{\ensuremath{\mathop{}\mathopen{}\mathcal{O}\mathopen{}(#1)}}
\newcommand{\bigOfrac}[1]{\ensuremath{\mathop{}\mathopen{}\mathcal{O}\mathopen{}\big(#1\big)}}
\newcommand{\bigOfracEq}[1]{\ensuremath{\mathop{}\mathopen{}\mathcal{O}\mathopen{}\Big(#1\Big)}}

\newcommand{\overbar}[1]{\mkern 1.5mu\overline{\mkern-1.5mu#1\mkern-1.5mu}\mkern 1.5mu}
\newcommand{\ud}{\mathrm{d}}
\newcommand{\RNum}[1]{\uppercase\expandafter{\romannumeral #1\relax}}
\newcommand{\numcircled}[1]{\tikz[baseline=(char.base)]{
            \node[shape=circle,draw,inner sep=0.4pt] (char) {#1};}}
\DeclareMathOperator*{\argmin}{arg\,min}

\newcolumntype{Y}{>{\centering\arraybackslash}X}

\makeatletter
\def\raisedotfill{%
  \leavevmode
  \cleaders \hb@xt@ .44em{\hss\raise0.5ex\hbox{.}\hss}\hfill
  \kern\z@}
\def\TPT@doparanotes{\par
   \prevdepth\z@ \TPT@hsize
   \TPTnoteSettings
   \parindent\z@ \pretolerance 8
   \linepenalty 200
   \renewcommand\item[1][]{\relax\ifhmode \begingroup
       \unskip
       \advance\hsize 10em 
       \penalty -45 \hskip\z@\@plus\hsize \penalty-19
       \hskip .15\hsize \penalty 9999 \hskip-.15\hsize
       \hskip .01\hsize\@plus-\hsize\@minus.01\hsize 
       \hskip .3em\@plus .3em
      \endgroup\fi
      \tnote{##1}\,\ignorespaces}%
   \let\TPToverlap\relax
   \def\endtablenotes{\par}%
}
\pretocmd\@bibitem{\color{black}\csname keycolor#1\endcsname}{}{\fail}
\newcommand\citecolor[1]{\@namedef{keycolor#1}{\color{blue}}}
\makeatother

\begin{document}
    \title{Distributed Nonconvex Optimization: Gradient-free Iterations and $\epsilon$-Globally Optimal Solution}

    \author{Zhiyu He\textsuperscript{$\dagger$,$\ddagger$}~\IEEEmembership{Student Member, IEEE}, Jianping He\textsuperscript{$\dagger$}~\IEEEmembership{Senior Member, IEEE}, Cailian Chen\textsuperscript{$\dagger$}~\IEEEmembership{Member, IEEE}, \\ and Xinping Guan\textsuperscript{$\dagger$}~\IEEEmembership{Fellow, IEEE} 
        \thanks{This work was supported in part by the National Natural Science Foundation of China under Grant 62373247, Grant 92167205, Grant 61933009, and Grant 62025305, and in part by the Max Planck ETH Center for Learning Systems. This paper was presented in part at the 2020 American Control Conference\cite{he2020cpca}. \emph{(Corresponding author: Jianping He.)} \\

        \textsuperscript{$\dagger$}Department of Automation, Shanghai Jiao Tong University, Shanghai 200240, China, Key Laboratory of System Control and Information Processing, Ministry of Education of China, Shanghai 200240, China, and Shanghai Engineering Research Center of Intelligent Control and Management, Shanghai 200240, China (email: hzy970920@sjtu.edu.cn; jphe@sjtu.edu.cn; cailianchen@sjtu.edu.cn; xpguan@sjtu.edu.cn). \\

        \textsuperscript{$\ddagger$}Automatic Control Laboratory, ETH Zurich, 8092 Zurich, Switzerland (email: zhiyhe@ethz.ch).
        }}

    \maketitle

    \begin{abstract}
        Distributed optimization utilizes local computation and communication to realize a global aim of optimizing the sum of local objective functions. This article addresses a class of constrained distributed nonconvex optimization problems involving univariate objectives, aiming to achieve global optimization without requiring local evaluations of gradients at every iteration. We propose a novel algorithm named CPCA, exploiting the notion of combining Chebyshev polynomial approximation, average consensus, and polynomial optimization. The proposed algorithm is i) able to obtain $\epsilon$-globally optimal solutions for any arbitrarily small given accuracy $\epsilon$, ii) efficient in both zeroth-order queries (i.e., evaluations of function values) and inter-agent communication, and iii) distributed terminable when the specified precision requirement is met. The key insight is to use polynomial approximations to substitute for general local objectives, distribute these approximations via average consensus, and solve an easier approximate version of the original problem. Due to the nice analytic properties of polynomials, this approximation not only facilitates efficient global optimization, but also allows the design of gradient-free iterations to reduce cumulative costs of queries and achieve geometric convergence for solving nonconvex problems. We provide a comprehensive analysis of the accuracy and complexities of the proposed algorithm.
    \end{abstract}

    
\section{Introduction}\label{sec:intro}
The developments of distributed optimization algorithms are motivated by wide application scenarios, including distributed learning, statistics, estimation, and coordination in large-scale networked systems. These algorithms enable agents in a network to collaboratively optimize a global objective function, which is generally the sum or average of local objectives, through local computations and communication. Owing to their features of exploiting network-wide resources and not requiring central coordinators, these algorithms enjoy higher efficiency, scalability, and robustness compared with their centralized counterparts\cite{nedic2018network}.

To solve coupled optimization problems by leveraging local interactions, the studies of distributed optimization benefit from the extensive research on optimization algorithms and consensus. Conversely, they also provide new insights, including the dependence of convergence rates on network topology\cite{nedic2018network} and the acceleration via gradient tracking\cite{shi2015extra,di2016next,xu2017convergence,qu2018harnessing}. There is a close relationship between the studies of distributed optimization and those of optimization and consensus.

\subsection{Related Work}
\textit{Distributed Convex Optimization}
algorithms generally fall into two categories, i.e., primal and dual-based methods. Primal methods combine gradient descent with consensus, thus driving local estimates to converge to the globally optimal point. Early methods exhibit sublinear convergence for nonsmooth convex objectives\cite{nedic2009subgradient}. By exploiting current and historical local gradients to track the global gradient, recent works achieve linear convergence rates for strongly convex and smooth objectives\cite{shi2015extra,nedic2017achieving,xu2017convergence,qu2018harnessing}. Current foci include stochastic gradients\cite{xin2021improved}, 
asynchronous computations \cite{pu2020push}, and second-order methods\cite{houska2016augmented}. 
Dual-based methods introduce consensus equality constraints\cite{makhdoumi2017convergence}, and then solve the dual problem\cite{scaman2017optimal} or carry on primal-dual updates\cite{shi2014linear,makhdoumi2017convergence}. Their extension to handle time-varying directed graphs is achieved by relaxing consensus constraints into asymptotically accurate closeness requirements\cite{aybat2019distributed,khatana2023dc}.

\textit{Distributed Nonconvex Optimization:}
Several noticeable distributed algorithms have been proposed, e.g., \cite{bianchi2012convergence,di2016next,tatarenko2017non,hong2017prox,wai2017decentralized,scutari2019distributed}. Their algorithmic frameworks share similarities with those for convex problems. Nevertheless, various techniques, including stochastic gradient descent\cite{bianchi2012convergence}, perturbations\cite{tatarenko2017non}, proximal methods\cite{hong2017prox}, polynomial filtering\cite{sun2019distributed}, and successive convex approximation\cite{di2016next,scutari2019distributed} enable agents to iteratively converge to stationary or locally optimal points of nonconvex problems. 

\textit{Zeroth-order Distributed Optimization}
is motivated by the concern that black-box procedures or resource limitations may inhibit access to the gradients of objective functions. To address this issue, the key idea is to utilize zeroth-order information (i.e., function evaluations) to construct randomized gradient estimates. Distributed zeroth-order algorithms \cite{hajinezhad2019zone,tang2019distributed} perform iterative updates based on these estimates. Their convergence rates match those of their first-order counterparts. 

\textit{Distributed Constrained Optimization:}
A common strategy for tackling convex local constraint sets is to project the newly generated solutions to the local set\cite{nedic2010constrained}. Further issues have been studied, including random projections and asynchronous updates. 
In other algorithms, the feasibility of new local estimates is satisfied by minimizing local surrogate functions\cite{di2016next,scutari2019distributed} or solving proximal minimization problems\cite{margellos2017distributed} right over the constraint sets. For general problems with inequality and equality constraints, the common practice is to introduce multipliers and then use consensus-based primal-dual gradient methods \cite{yi2020distributed}. 
Further, a multitude of second-order algorithms feature local superlinear convergence rates and handle nonconvex objectives and constraints. The key idea is to locally solve a sequence of nonlinear optimization problems and achieve coordination by solving a convex quadratic program. This coordination problem can be addressed in a centralized (as in ALADIN \cite{houska2016augmented}) or decentralized fashion\cite{engelmann2020decomposition,engelmann2021essentially}. 



\subsection{Motivations}
The above distributed optimization algorithms are provably efficient.
Nonetheless, there remain two notable unresolved issues. First, for problems with nonconvex objectives, only the convergence to stationary points is guaranteed. The limit point that these descent methods converge to is where the gradient of the global objective vanishes. It is hard to determine if this point is globally optimal, locally optimal, or a saddle point. Second, the load of oracle queries (i.e., evaluations of gradients or function values) generally grows with the number of iterations. This increasing load stems from the iterative algorithmic structure, where local oracle queries are constantly performed at every iteration. When the numbers of iterations are large (e.g., optimization over large-scale networks) and such queries are costly (e.g., in hyperparameter optimization\cite{yang2020hyperparameter}), 
the cumulative load of queries can be a critical bottleneck. To resolve these issues, we need to find a new path different from those of existing distributed algorithms.

We are inspired by the close link between function approximation and optimization. Various studies introduce approximation to improve the performance of optimization algorithms, including speeding up convergence (e.g., Newton's method\cite{boyd2004convex}) and reducing computational costs (e.g., successive convex approximation techniques\cite{di2016next,scutari2019distributed}). 
These algorithms share a common feature of setting new estimates as minimizers of some local approximations of the objective function.
The minimizers of these approximations are readily available in general, and they can often be expressed in closed forms. 
The related distributed implementations include \cite{varagnolo2015newton,di2016next,scutari2019distributed}. 

To exploit global approximations, some work on numerical analysis uses Chebyshev polynomial approximations to substitute for univariate functions defined on intervals\cite{trefethen2013approximation,driscoll2014chebfun,boyd2014solving}. These approximations are polynomial interpolations in Chebyshev points within intervals. This substitution contributes to studying the properties of functions, including optimization, root-finding, and integration. For a general nonconvex objective, as long as it is Lipschitz continuous, we can construct its arbitrarily precise approximation (i.e., proxy) on the entire interval. Therefore, we turn to solve an easier problem of optimizing this proxy for the global objective. Moreover, the vector of coefficients of a local proxy serves as its compact representation. We can then distribute these vectors via average consensus to let agents obtain the needed global proxy.

The above observations motivate the new algorithm presented in this article.
Instead of constantly querying first-order or zeroth-order information at every iteration, we utilize function evaluations (i.e., zeroth-order information) to construct polynomial approximations for local objectives in advance, thus facilitating subsequent information dissemination and optimization. This design helps to achieve efficient optimization of the global nonconvex objective function, with lower costs in zeroth-order queries and inter-agent communication.

\subsection{Contributions}
We propose a Chebyshev-Proxy-and-Consensus-based Algorithm (CPCA) to solve optimization problems with nonconvex univariate objectives and convex constraint sets in a distributed manner.
The main contributions are as follows.


\begin{itemize}
\item We develop a novel algorithm CPCA featuring gradient-free, consensus-based iterations. The key idea is to use polynomial approximations to substitute for general objectives, distribute these approximations via average consensus, and solve instead an approximate version of the original problem. It differs from existing iterative distributed methods and offers a new perspective on addressing distributed optimization problems.
 
\item We prove that CPCA obtains the $\epsilon$-globally optimal solution of nonconvex problems, where $\epsilon$ is any arbitrarily small given solution accuracy.  
To this end, we provide a general rule, which tightly connects and controls the error bounds for the inner stages of the proposed algorithm. Furthermore, we show that it achieves distributed termination once the precision requirement is met.

\item We characterize the complexities of zeroth-order queries, floating-point operations (i.e., \textit{flops}\footnote{A flop is defined as one addition, subtraction, multiplication, or division of two floating-point numbers\cite{boyd2004convex}.}), and inter-agent communication of the proposed algorithm. Thanks to its unique introduction of approximation and gradient-free iterations, CPCA is efficient in terms of communication rounds and queries, see the comparison in Section~\ref{subsec:complexity}.
\end{itemize}

The differences between this article and its conference version\cite{he2020cpca} include \romannumeral1) the incorporation of the distributed stopping mechanism for consensus iterations, \romannumeral2) another reformulation based on semidefinite programming to optimize the polynomial proxy, and \romannumeral3) results and discussions on the multivariate extension. 

\subsection{Organization and Notations}
The rest of this article is organized as follows. Section~\ref{sec:formulation} describes the problem of interest and provides preliminaries. Section~\ref{sec:algorithm} develops the algorithm CPCA\@. Section~\ref{sec:analysis} analyzes the accuracy and complexities of the proposed algorithm. Further discussions on application scenarios and the algorithmic structure are given in Section~\ref{sec:discussion}. Section~\ref{sec:experiment} presents the simulation results. Finally, Section~\ref{sec:conclusion} concludes this article. 

\emph{Notations:}
let $\lVert a \rVert$ and $\lVert a \rVert_{\infty}$ be the $\ell_{2}$-norm and $\ell_{\infty}$-norm of a vector $a\in \mathbb{R}^{n}$, respectively. 
The superscript $t$ denotes the number of iterations. The subscripts $i,j$ denote the indexes of agents. The script $k$ in parentheses denotes the index of elements in a vector. 

\section{Problem Description and Preliminaries}\label{sec:formulation}

\subsection{Problem Description}
Consider a network with $N$ agents, each of which owns a local objective function $f_i(x):X_{i} \to \mathbb{R}$ and a local constraint set $X_{i} \subset \mathbb{R}$. The goal is to solve the following constrained optimization problem
\begin{equation}\label{problem:main_focus}
    \begin{split}
        \min_{x} \quad & f(x) = \frac{1}{N}\sum_{i=1}^{N} f_i(x) \\
        \textrm{s.t.} \quad & x \in X = \bigcap_{i=1}^{N} X_{i}
    \end{split}
\end{equation}
in a distributed manner. The network is described by an undirected connected graph $\mathcal{G}=(\mathcal{V},\mathcal{E})$, where $\mathcal{V}$ is the set of agents and $\mathcal{E} \subseteq \mathcal{V} \times \mathcal{V}$ is the set of edges. Agent $j$ can receive information from agent $i$ if and only if (\textit{iff}) $(i,j) \in \mathcal{E}$. Two basic assumptions are given as follows.

\begin{assumption}\label{assump:lipschitz_continous}
    The local objective $f_{i}(x)$ is Lipschitz continuous on $X_{i}$.
\end{assumption}


\begin{assumption}\label{assump:constraint_set}
    The local constraint set $X_{i}$ is a closed, bounded, and convex set.
\end{assumption}

Assumptions \ref{assump:lipschitz_continous} and \ref{assump:constraint_set} are general, commonly seen in the literature, and satisfied by typical problems of practical interests (see \cite{nedic2010constrained,nedic2018network,scutari2019distributed} and the references therein).

Problem \eqref{problem:main_focus} involves nonconvex objectives and convex constraint sets, and, therefore, it is a constrained nonconvex optimization problem. Under Assumption \ref{assump:constraint_set}, for all $i \in \mathcal{V}$, $X_{i}$ is a closed interval. Thus, let $X_{i} = [a_{i},b_{i}]$, where $a_{i},b_{i} \in \mathbb{R}$. It follows that $X=[a,b]$, where $a=\max_{i \in \mathcal{V}} a_{i},~b=\min_{i \in \mathcal{V}} b_{i}$.


General distributed optimization algorithms with gradient tracking\cite{shi2015extra,di2016next,xu2017convergence,qu2018harnessing} can be abstracted by
\begin{subequations}\label{eq:generalDGT}
\begin{align}
    x_{i}^{t+1} &= \mathcal{F}_t\left({\textstyle \sum_j} w_{ij} x_{j}^{t}, s_{i}^{t}\right), \\
    s_{i}^{t+1} &= {\textstyle \sum_j} w_{ij}' s_{j}^{t} + \nabla f_{i}(x_{i}^{t+1}) - \nabla f_{i}(x_{i}^{t}),
\end{align}
\end{subequations}
where $\mathcal{F}_t$ is the update rule, and $w_{ij}, w_{ij}' \geq 0$ are the weights\cite{qu2018harnessing}. Note that for every agent $i$, the evaluation of the local gradient $\nabla f_i$ is performed at every iteration. As discussed in Section~\ref{sec:intro}, even for univariate problems, the cumulative cost of such evaluations can be high when certain problems (e.g., hyperparameter optimization) are solved over large-scale networks. Additionally, for nonconvex univariate problems, algorithms as \eqref{eq:generalDGT} only ensure convergence to stationary points\cite{scutari2019distributed,tang2019distributed}. In contrast, in this article we aim to overcome the aforementioned restrictions by taking a new path of introducing polynomial approximation rather than following the existing framework \eqref{eq:generalDGT}.

\subsection{Preliminaries}

\emph{Consensus algorithms} enable agents in a networked system to reach global agreement via local information exchange only.
Let $\mathcal{N}_{i} = \{j|(j,i)\in \mathcal{E}\}$ be the set of neighbors of agent $i$ and $\operatorname{deg}(i) = |\mathcal{N}_{i}|$ be its degree, where $|\mathcal{N}_{i}|$ denotes the cardinality of $\mathcal{N}_{i}$. Suppose that every agent $i$ maintains a local variable $x_{i}^{t}\in \mathbb{R}$. 
The maximum consensus algorithm\cite{saber2003consensus} is
\begin{equation}\label{eq:max_consensus}
    x_i^{t+1} = \max_{j \in \mathcal{N}_{i}} x_j^{t}.
\end{equation}
It has been proven that with \eqref{eq:max_consensus}, all $x_i^{t}$ converge exactly to $\max_{i\in \mathcal{V}} x_{i}^{0}$ in $T~(T \leq D)$ iterations, where $D$ is the diameter (i.e., the greatest distance between any pair of agents) of $\mathcal{G}$. 
The average consensus algorithm based on lazy Metropolis weights\cite{nedic2018network} is
\begin{equation}\label{eq:avg_consensus}
    x_i^{t+1} = x_i^{t}+\frac{1}{2} \sum_{j\in \mathcal{N}_i} \frac{x_{j}^{t}-x_{i}^{t}}{\max (\operatorname{deg}(i),\operatorname{deg}(j))}.
\end{equation}
This algorithm requires an agent to exchange its local variable $x_{i}^{t}$ and degree $\operatorname{deg}(i)$ with its neighbors. With \eqref{eq:avg_consensus}, all $x_{i}^{t}$ converge geometrically to the average $\overbar{x} = 1/N \sum_{i=1}^{N} x_{i}^{0}$ of the initial values\cite{nedic2018network}. 
Let $x^{t} = [x_{1}^{t},\ldots,x_{N}^{t}]^{\top}$, and let $t(\xi)$ denote the first time $t$ when $\frac{\|x^{t} - \overbar{x}1\|}{\|x^{0} - \overbar{x}1\|} \leq \xi$. Then, by referring to \cite[Proposition 5]{nedic2018network}, we have
\begin{equation}\label{eq:convergence_time}
    t(\xi) \sim \bigOfracEq{\log \frac{1}{\xi}}.
\end{equation}

\emph{Chebyshev polynomial approximation} uses truncated Chebyshev series to approximate functions, thus facilitating numerical analysis. 
For a Lipschitz continuous function $g(x)$ on $[a,b]$, its degree $m$ Chebyshev approximation $p^{(m)}(x)$ is 
\begin{equation}\label{eq:cheb_rep}
    p^{(m)}(x) = \sum_{j=0}^{m} c_{j} T_{j} \left(\frac{2x-(a+b)}{b-a}\right), \quad x\in [a,b],
\end{equation}
where $c_{j} \in \mathbb{R}$ is the Chebyshev coefficient, and $T_{j}(\cdot)$ is the $j$-th Chebyshev polynomial on $[-1,1]$ that satisfies $\forall u\in [-1,1], \forall j, |T_{j}(u)|\leq 1$ and forms a basis\cite{trefethen2013approximation}. As $m$ increases, $p^{(m)}(x)$ converges to $g(x)$ uniformly on the given interval\cite{trefethen2013approximation}. 
The convergence rates depend on the smoothness of $g(x)$ and are characterized in Section~\ref{subsec:DegreeApprox}. The above observation makes computing $p^{(m)}(x)$ a practical way to construct an arbitrarily close polynomial approximation of $g(x)$, as theoretically guaranteed by the \emph{Weierstrass Approximation Theorem}\cite[Theorem 6.1]{trefethen2013approximation}.

\emph{Sum of Squares Polynomials:}
Let the vector of Chebyshev polynomials of degrees until $d$ be $v(x)\triangleq [T_{0}(x),\ldots,T_{d}(x)]^{\top}$. The elements of $v(x)$ constitute an orthogonal basis for all polynomials of degree $d$ or less. Any univariate polynomial $p(x)$ of degree $2d$ is non-negative (or equivalently, a sum of squares (SOS)) on $\mathbb{R}$ \textit{iff} there exists a positive semidefinite matrix $Q$ of order $d+1$ (i.e., $Q\in \mathbb{S}^{d+1}_{+}$) such that
\begin{equation}\label{eq:sos_psd}
    p(x) = v(x)^{\top}Qv(x) = \sum_{i,j=0}^{d} Q_{ij} T_{i}(x)T_{j}(x),
    \end{equation}
where the rows and columns of $Q$ are indexed by $0,1,\ldots,d$.
We refer readers to \cite[Lemma 3.33]{blekherman2012semidefinite} for details. 
\section{CPCA: Design and Analysis}\label{sec:algorithm}
\begin{figure}[tb]
\begin{center}
    \includegraphics[width=\columnwidth]{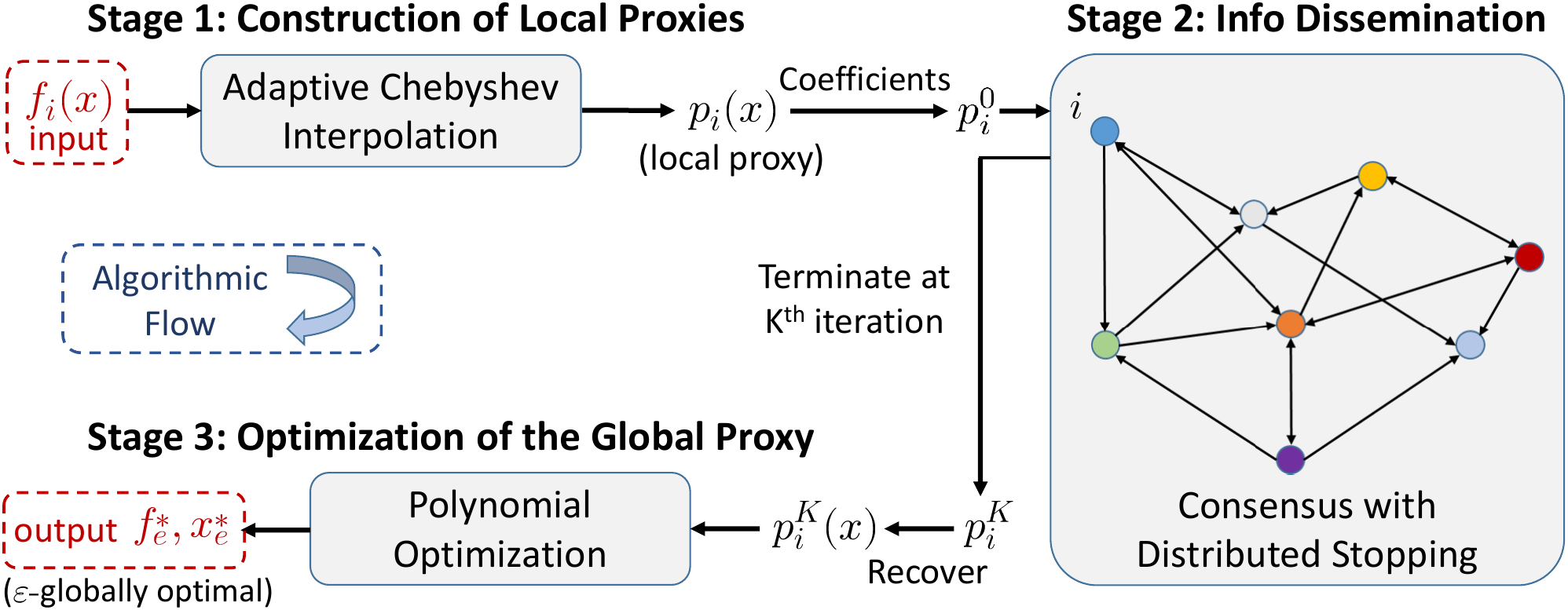}
    \caption{An overview of CPCA and its algorithmic flow.}
    \label{fig:overview}
\end{center}
\end{figure}

We present CPCA to solve problem \eqref{problem:main_focus} in a distributed manner. Fig.~\ref{fig:overview} illustrates the main architecture and information flow of the proposed algorithm.

\subsection{Construction of Local Chebyshev Proxies}\label{subsec:init}
In this stage, every agent $i$ constructs a polynomial approximation $p_{i}(x)$ of $f_{i}(x)$ on $X=[a,b]$ to achieve
\begin{equation}\label{eq:approx_require}
    \mathbf{G_1}: \quad |f_{i}(x) - p_{i}(x)| \leq \epsilon_{1}, \quad \forall x\in [a,b],
\end{equation}
where $\epsilon_{1} > 0$ is a specified error bound. Specifically, Chebyshev polynomial approximation is chosen because \romannumeral1) it is in the global sense, which contrasts with the Taylor expansion that characterizes the local landscape; \romannumeral2) it allows compact representations of complex functions via vectors of coefficients; and \romannumeral3) the ranges of basis polynomials are identical and bounded, which facilitates establishing the accuracy guarantee.
The key insight is to adaptively settle at a suitable degree of the Chebyshev approximation when a stopping criterion is satisfied\cite{boyd2014solving}.
Agent $i$ initializes $m_{i} = 2$ and calculates a Chebyshev interpolant of degree $m_{i}$. It first evaluates $f_{i}(x)$ at the set $S_{m_{i}} \triangleq \{x_{0},\ldots,x_{m_{i}}\}$ of $\left(m_{i}+1\right)$ points by
\begin{equation}\label{eq:grid_point_evaluation}
        x_{k} = \frac{b-a}{2}\cos\left(\frac{k\pi}{m_{i}} \right)+ \frac{a+b}{2}, \quad f_{k} = f_{i}(x_{k}),
\end{equation}
where $k=0,1,\dots,m_{i}$. Then, based on these function evaluations, it computes the Chebyshev coefficients of the interpolant of degree $m_{i}$ by
\begin{equation}\label{eq:coeff_value_relation}
    c_{j} = \frac{1}{m_{i}} \left(f_{0} \! + \! f_{m_{i}}\cos(j\pi)\right) + \frac{2}{m_{i}} \sum_{k=1}^{m_{i}-1} f_{k}\cos\left(\frac{jk\pi}{m_{i}}\right),
\end{equation}
where $j=0,1,\dots,m_{i}$. 
The degree $m_{i}$ is doubled at every iteration until the stopping criterion
\begin{equation}\label{eq:adapt_stop_rule}
    \max_{x_{k}\in \left(S_{2m_{i}}-S_{m_{i}}\right)} |f_{i}(x_{k})-p_{i}(x_{k})| \leq \epsilon_{1}
\end{equation}
is met, where $p_{i}(x)$ is in the form of \eqref{eq:cheb_rep} with $\{c_{j}\}$ being the coefficients.
To check \eqref{eq:adapt_stop_rule}, we sample points in $S_{2m_i} - S_{m_i}$ (i.e., the points that are in $S_{2m_i}$ but are not in $S_{m_i}$).
Since $S_{m_{i}}\subset S_{2m_{i}}$, the evaluations of $f_{i}(x)$ are continuously reused. Note that agents know the intersection set $X=[a,b]$ by running a finite number of max/min consensus iterations as \eqref{eq:max_consensus} in advance. 
The aforementioned procedures return $p_{i}(x)$ satisfying (8) for most Lipschitz continuous $f_{i}(x)$ encountered in practice\cite{boyd2014solving}. To handle special cases where (8) may fail, in numerical practice, we can additionally use the chopping technique\cite{aurentz2017chopping} that examines the plateau of Chebyshev coefficients near machine precision ($\sim 10^{-16}$) and obtains a satisfactory approximation through appropriate chopping.

\subsection{Consensus-based Information Dissemination}\label{subsec:iterations}
In this stage, agents perform consensus-based iterations to update their local variables $p_{i}^{t}$, whose initial value $p_{i}^{0} = [c_{0},\ldots,c_{m_{i}}]^{T}$ stores the coefficients of the local approximation $p_{i}(x)$. The goal is to enable $p_{i}^{t}$ to converge to the average of the initial values $\overbar{p}= \frac{1}{N} \sum_{i=1}^{N} p_{i}^{0}$, since $\bar{p}$ leads to a satisfactory approximation of the global objective $f(x)$ in \eqref{problem:main_focus}. Specifically, let $K$ be the total number of iterations and $p_{i}^{K}$ be the local variable of agent $i$ at time $K$. We aim to ensure that
\begin{equation*}
    \mathbf{G_2}: \quad \max_{i \in \mathcal{V}} \left \lVert p_{i}^{K} -\overbar{p} \right \rVert _{\infty} \leq \delta \triangleq \frac{\epsilon_{2}}{m+1}, \quad m \triangleq \max_{i \in \mathcal{V}} m_{i}
\end{equation*}
holds, where $\epsilon_{2} > 0$ is an error bound and $m$
is the highest of the degrees of local approximations. Then, a sufficiently precise proxy (recovered from $p_{i}^{K}$) for the global objective is acquired. The challenge is fulfilling in-time termination when the accuracy requirement is met. The details are as follows.
 
Every agent updates its local variable $p_{i}^{t}$ according to 
\begin{equation}\label{eq:avg_consensus_implement}
    p_i^{t+1} = p_i^{t}+\frac{1}{2} \sum_{j\in \mathcal{N}_i} \frac{p_{j}^{t}-p_{i}^{t}}{\max (\operatorname{deg}(i),\operatorname{deg}(j))}.
\end{equation}
To achieve distributed stopping once the specified precision is reached, we incorporate the max/min-consensus-based stopping mechanism in \cite{yadav2007distributed}. The main idea is to exploit the monotonicity of the sequences of maximum and minimum local variables across the network and the finite-time convergence of max/min consensus algorithms. This mechanism requires some knowledge of the diameter $D$ of $\mathcal{G}$.
\begin{assumption}\label{assump:upperbound}
    Every agent knows an upper bound $U$ on $D$.
\end{assumption}

The distributed estimation of $D$ can be realized via the Extrema Propagation technique\cite{jesus2014survey}. If Assumption \ref{assump:upperbound} holds, agents will know how long to wait to ensure that max/min consensus algorithms converge.

To realize distributed stopping, there are two auxiliary variables $r_{i}^{t}$ and $s_{i}^{t}$, which are initialized as $p_{i}^{0}$. These variables are updated in parallel with $p_{i}^{t}$ by
\begin{equation}\label{eq:max_consensus_auxiliary}
    r_{i}^{t+1}(k) = \max_{j\in \mathcal{N}_{i}} r_{j}^{t}(k), \quad s_{i}^{t+1}(k) = \min_{j\in \mathcal{N}_{i}} s_{j}^{t}(k),
\end{equation}
where $k=1,\ldots,m+1$, and they are reinitialized as $p_{i}^{t}$ every $U$ iterations to facilitate the constant dissemination of recent information on $p_{i}^{t}$. At time $t_{l}=lU(l\in \mathbb{N})$, every agent $i$ periodically checks the following criterion
\begin{equation}\label{eq:itr_stop_rule}
    \|r_{i}^{t_{l}} - s_{i}^{t_{l}}\|_{\infty} \leq \delta.
\end{equation}
If \eqref{eq:itr_stop_rule} is met at time $K$, then agents terminate the iterations. Consequently, $p_{i}^{K}$ is sufficiently close to $\overbar{p}$, as characterized by the following theorem.
\begin{theorem}\label{thm:itr_precision}
    When \eqref{eq:itr_stop_rule} is satisfied, we have 
    \begin{equation}\label{eq:itr_precision_requirement}
        \max_{i \in \mathcal{V}} \left \lVert p_{i}^{K} -\overbar{p} \right \rVert _{\infty} \leq \delta, \quad \delta = \frac{\epsilon_{2}}{m+1}.
    \end{equation}
\end{theorem}
\begin{proof}
    See Appendix~\ref{subsec:proof_itr}. 
\end{proof}

\begin{remark}\label{rem:agree_dim}
    While performing updates by \eqref{eq:avg_consensus} and \eqref{eq:max_consensus_auxiliary}, agents ensure agreements in dimension by aligning and padding zeros to those local variables of smaller sizes. After a finite number of iterations (less than $D$), all these local variables will be of the same size $\max_{i\in \mathcal{V}} m_{i}(= m) + 1$. Afterward, the exchanged number of coefficients is fixed.
\end{remark}

\begin{remark}
    The adaptive approximation in Section~\ref{subsec:init} offers a compact way of representation and exchange. In comparison, the strategy of directly exchanging local objectives can be hard to implement due to the challenge of encoding and transmitting heterogeneous objectives. Moreover, the method of using a fixed, high degree may be inefficient when the accuracy requirement is medium and a moderate degree suffices.
\end{remark}

\subsection{Polynomial Optimization via Finding Stationary Points}\label{subsec:alternative_poly_opt}
In this stage, agents optimize the polynomial proxy $p_{i}^{K}(x)$ recovered from $p_{i}^{K}$ independently, i.e.,
 \begin{equation*}
     \mathbf{G_3}: \quad \min_x ~ p_i^K(x) \quad \textrm{s.t.}~x \in [a,b],
 \end{equation*}
 thus obtaining $\epsilon$-optimal solutions to problem \eqref{problem:main_focus}. To achieve this goal, we first obtain all the stationary points of $p_{i}^{K}(x)$ by calculating the eigenvalues of a colleague matrix constructed from its Chebyshev coefficients\cite{trefethen2013approximation}. Then, by comparing the values of $p_{i}^{K}(x)$ at these critical points, we decide the optimal value and the set of optimal points of $p_{i}^{K}(x)$. 


Suppose that $p_{i}^{K}(x)$ is in the form of \eqref{eq:cheb_rep} with $\{c_{j}'|j=0, \ldots, m\}$ being the coefficients. Then,
\begin{equation*}
    \frac{\ud p_{i}^{K}(x)}{\ud x} = \sum_{j=0}^{m-1} \tilde{c}_{j} T_{j} \left(\frac{2x-(a+b)}{b-a}\right), \quad x\in (a,b).
\end{equation*}
The coefficients $\{\tilde{c}_{j}|j=0,\ldots,m-1\}$ are obtained from the following recurrence formula
\begin{equation}\label{eq:recurrence_coeff}
\tilde{c}_{j} = \left\{
    \begin{aligned}
        & \tilde{c}_{j+2}+2(j+1)Sc_{j+1}', & & j = m-1,\ldots,1, \\
        & {\textstyle \frac{1}{2}} \tilde{c}_{2}+S c_{1}', & & j = 0,
    \end{aligned}
    \right.
\end{equation}
where $S=2/(b-a), \tilde{c}_{m}=\tilde{c}_{m+1}=0$. By \cite{trefethen2013approximation}, the roots of $\ud p_{i}^{K}(x)/\ud x$ (i.e., the stationary points of $p_{i}^{K}(x)$) are the eigenvalues of the following square colleague matrix $M_{C}$ of order $m-1$
\begin{equation*}\label{eq:colleague_matrix}
    \begin{bmatrix}
    0 & 1 & & & \\
    \frac{1}{2} & 0 & \frac{1}{2} & & \\
    & \ddots & \ddots & \ddots & \\
    & & \frac{1}{2} & 0 & \frac{1}{2} \\
    -\frac{\tilde{c}_{0}}{2\tilde{c}_{m-1}} & -\frac{\tilde{c}_{1}}{2\tilde{c}_{m-1}} &  \cdots & \frac{1}{2}-\frac{\tilde{c}_{m-3}}{2\tilde{c}_{m-1}} & -\frac{\tilde{c}_{m-2}}{2\tilde{c}_{m-1}}
    \end{bmatrix}
\end{equation*}
Let $E$ be the set of all the real eigenvalues of $M_{C}$ that lie in $X=[a,b]$. Then, the optimal value $f_{e}^{*}$ and the set of optimal points $X_{e}^{*}$ of $p_{i}^{K}(x)$ are
\begin{equation}\label{eq:opt_selection}
    f_{e}^{*} = \min_{x \in X_{K}} p_{i}^{K}(x), \quad
    X_{e}^{*} = \argmin_{x \in X_{K}} p_{i}^{K}(x), 
\end{equation}
where $X_{K} = E \cup \{a,b\}$ is the set of all the critical points.



\subsection{Description of CPCA}\label{subsec:algorithmic_description}
CPCA consists of the above three stages and is summarized as Algorithm~\ref{alg:alg_overall}, where lines \ref{procedure:init}-\ref{procedure:init_end} perform the initialization; lines \ref{procedure:calc_lc_prxy}-\ref{procedure:calc_lc_prxy_end} construct the local proxy; lines \ref{procedure:average_consensus}-\ref{procedure:average_consensus_end} are consensus iterations with distributed stopping; and line \ref{procedure:poly_opt_slv} is polynomial optimization.
With CPCA, we achieve goals $\mathbf{G_1}$, $\mathbf{G_2}$, and $\mathbf{G_3}$, thereby solving problem~\eqref{problem:main_focus}. Furthermore, the algorithm leverages the given error tolerance $\epsilon > 0$ to set the error bounds (i.e., $\epsilon_{1}$ and $\epsilon_{2}$) utilized in the corresponding stages. To meet the specified accuracy requirement, these bounds satisfy
\begin{equation}\label{eq:sum_err_require}
    \epsilon_{1}+\epsilon_{2}=\epsilon, \qquad \epsilon_{1},\epsilon_{2}>0.
\end{equation}
We set these error bounds to be $\epsilon/2$ in Algorithm~\ref{alg:alg_overall}. 

\begin{algorithm}[tb]
\small
\caption{CPCA}
\label{alg:alg_overall}
    \begin{algorithmic}[1]
    \REQUIRE $f_{i}(x),X_{i}=[a_{i},b_{i}],U,\text{ and }\epsilon$.
    \ENSURE $f^{*}_{e}$ for every agent $i \in \mathcal{V}$.
    \STATE {\bfseries Initialize:} $a_{i}^{0}=a_{i},b_{i}^{0}=b_{i},m_{i}=2$. \label{procedure:init}
    \FOR{{\bfseries each} agent $i\in \mathcal{V}$}
        \FOR{$t=0, \ldots, U-1$}
            \STATE $\displaystyle a_{i}^{t+1} = \max_{j\in \mathcal{N}_{i}} a_{j}^{t}, ~ b_{i}^{t+1} = \min_{j\in \mathcal{N}_{i}} b_{j}^{t}$.
        \ENDFOR
        \STATE Set $a=a_{i}^{t},~b=b_{i}^{t}$. \label{procedure:init_end}
        \STATE Calculate $\{x_{j}\}\textrm{ and }\{f_{j}\}$ by \eqref{eq:grid_point_evaluation}. \label{procedure:interpolation} \label{procedure:calc_lc_prxy}
        \STATE Calculate $\{c_{k}\}$ by \eqref{eq:coeff_value_relation}.
        \STATE If \eqref{eq:adapt_stop_rule} is satisfied (where $\displaystyle \epsilon_{1} = \epsilon/2$), go to step \ref{procedure:construct_vector}. Otherwise, set $m_{i} \leftarrow 2m_{i}$ and go to step \ref{procedure:interpolation}. \label{procedure:calc_lc_prxy_end}
        \STATE Set $p_{i}^{0}=r_{i}^{0}=s_{i}^{0}= [c_{0},c_{1},\dots,c_{m_{i}}]^{\top},l=1$. \label{procedure:construct_vector} \label{procedure:average_consensus}
        \FOR{$t = 0,1,\dots$}
            \IF{$t=lU$}
                \IF{$l=1$}
                    \STATE $\delta = \epsilon_{2}/(m+1) = \epsilon/2(m+1)$.
                \ENDIF
                \IF{$\|r_{i}^{t} - s_{i}^{t}\|_{\infty} \leq \delta$}
                    \STATE $p_{i}^{K} = p_{i}^{t}$. ~\textbf{break}
                \ENDIF
                \STATE $r_{i}^{t}=s_{i}^{t}=p_{i}^{t},~l \leftarrow l+1$.
            \ENDIF
            \STATE Update $p_{i}^{t+1},r_{i}^{t+1},s_{i}^{t+1}$ by \eqref{eq:avg_consensus_implement} and \eqref{eq:max_consensus_auxiliary}.
            \STATE Set $t \leftarrow t+1$.
        \ENDFOR \label{procedure:average_consensus_end}
        \STATE Optimize $p_i^K(x)$ in \eqref{eq:opt_selection} and return $f^{*}_{e}, x_e^*$. \label{procedure:poly_opt_slv}
    \ENDFOR
    \end{algorithmic}
\end{algorithm}

The rule \eqref{eq:sum_err_require} allows the balance between computations and communication by adjusting the corresponding error bounds. For example, if the costs of inter-agent communication and intra-agent computations are high and low, respectively, then we can decrease the accuracy of iterations (i.e., increase $\epsilon_{2}$) and increase the precision of calculating approximations (i.e., decrease $\epsilon_{1}$). Consequently, less burden will be placed upon communication and more computational resources will be utilized, thus improving efficiency and adaptability.

The proposed algorithm differs from most distributed optimization algorithms (e.g., \cite{tang2019distributed,scutari2019distributed}) in both the shared information and the structure of iterations. First, existing algorithms generally exchange local estimates of optimal solutions, while CPCA communicates vectors of coefficients of local approximations (i.e., local proxies). Second, most existing algorithms call for local evaluations of gradients or function values at every iteration to update local estimates, whereas CPCA is free from such evaluations and only performs consensus-based updates related to vectors of coefficients. These updates lead to the acquisition of the global proxy and hence help to efficiently obtain $\epsilon$-globally optimal solutions, as we will see in Theorem~\ref{thm:alg_accuracy_result} in Section~\ref{sec:analysis}.

\section{Performance Analysis}\label{sec:analysis}

\subsection{Accuracy}
We first provide a lemma stating that if two functions $f,g: [a,b]\to \mathbb{R}$ are close on the entire interval, their optimal values $f(x_{f}^{*}) \text{ and } g(x_{g}^{*})$ are also close, where $x_f^*$ and $x_g^*$ are the optimal points of $f$ and $g$, respectively.
\begin{lemma}\label{lem:value_distance}
    If $f,g$ satisfy $|f(x)-g(x)|\leq \epsilon, \forall x\in [a,b]$, then
    \begin{equation}\label{eq:closeness}
        \big|f(x_{f}^{*})-g(x_{g}^{*})\big| \leq \epsilon.
    \end{equation}
    \begin{proof}
        We know that $|f(x_{f}^{*})-g(x_{f}^{*})| \leq \epsilon$. Since $x_{g}^{*}$ is the optimal point of $g$, we have $g(x_g^*)\leq g(x_f^*)$. Hence,
        \begin{equation*}
            g(x_g^*) \leq g(x_f^*) \leq f(x_f^*)+\epsilon,
        \end{equation*}
        which implies that $f(x_f^*)-g(x_g^*) \geq -\epsilon$. Similarly, we have
        \begin{equation*} 
            f(x_f^*) \leq f(x_g^*) \leq g(x_g^*)+\epsilon,
        \end{equation*}
        which leads to $f(x_f^*)-g(x_g^*) \leq \epsilon$. Therefore, \eqref{eq:closeness} holds.
    \end{proof}
\end{lemma}

Now we establish the accuracy of CPCA\@. The key idea is to prove the closeness between the global objective $f(x)$ in \eqref{problem:main_focus} and the locally recovered proxy $p_i^K(x)$. We use $\epsilon >0 \textrm{ and } f^{*}$ to denote the specified solution accuracy and the optimal value of problem~\eqref{problem:main_focus}, respectively.

\begin{theorem}\label{thm:alg_accuracy_result}
Suppose that Assumptions \ref{assump:lipschitz_continous}-\ref{assump:upperbound} hold. If $\epsilon > 0$ is specified, CPCA ensures that every agent obtains an $\epsilon$-optimal solution $f_{e}^{*}$ to problem \eqref{problem:main_focus}, i.e., $\left|f_{e}^{*}-f^{*}\right| \leq \epsilon$.
\begin{proof}
    See Appendix~\ref{subsec:proof_accuracy}. 
\end{proof}
\end{theorem}



\subsection{Complexity}\label{subsec:complexity}
We analyze the computational and communication complexities of the proposed algorithm. The following theorem establishes the orders of the numbers of zeroth-order queries (i.e., evaluations of values of local objective functions), inter-agent communication, and floating-point operations (i.e., flops) required for one agent. Suppose that one evaluation of the function value costs $F_{0}$ flops. In practice, this cost depends on the specific forms of objective functions\cite{boyd2004convex}.

\begin{theorem}\label{thm:alg_whole_complx}
    With CPCA, every agent obtains an $\epsilon$-optimal solution to problem~\eqref{problem:main_focus}, with $\bigO{m}$ zeroth-order queries, $\bigO{\log \frac{m}{\epsilon}}$ rounds of inter-agent communication, and $\bigOfrac{m\cdot \max(m, \log \frac{m}{\epsilon}, F_{0})}$ flops.
\end{theorem}
\begin{proof}
    The proof is provided in Appendix~\ref{subsec:proof_complexity}. We summarize the costs of zeroth-order queries, inter-agent communication, and flops of the proposed algorithm in Table~\ref{table:complex_result}. 
\end{proof}

\begin{table}[!tb]
\centering
\renewcommand \arraystretch{1.1}
\caption{Complexities of CPCA}
\label{table:complex_result}
\begin{threeparttable}
    \begin{tabularx}{\columnwidth}{*{3}{Y} c}
        \toprule
        {\bfseries Stages} & {\bfseries Queries} & {\bfseries Comm.} & {\bfseries Flops} \\
        \midrule
        init & $\bigO{m}$\tnote{1} & $/$ & $\bigO{m\cdot\max(m,F_{0})}$ \\
        iteration & $/$ & $\bigOfrac{\log \frac{m}{\epsilon}}$ & $\bigOfrac{m \log \frac{m}{\epsilon}}$ \\
        solve & $/$ & $/$ & $\bigO{m^2}$ \\
        \midrule
        whole & $\bigO{m}$ & $\bigOfrac{\log \frac{m}{\epsilon}}$ & $\bigOfrac{m\cdot \max(m, \log \frac{m}{\epsilon}, F_{0})}$ \\
        \bottomrule
        \addlinespace[0.5ex]
    \end{tabularx}
    \begin{tablenotes}[para]\footnotesize
        \item[1] The dependence of $m$ on $\epsilon$ is discussed in Lemma~\ref{lem:degree_approx}, Section~\ref{subsec:DegreeApprox}.
    \end{tablenotes}
\end{threeparttable}
\end{table}

Table~\ref{table:compare_result} compares the complexities of different algorithms. Though the algorithms in \cite{scutari2019distributed} and \cite{tang2019distributed} handle multivariate problems, the comparison here is for the univariate case. In CPCA, the evaluations of function values are not required during iterations and are only performed when approximations are constructed. This design implies that the number of evaluations will not increase with the number of iterations. Hence, the cumulative costs of queries can be significantly reduced especially for large-scale networks. Due to the exchange of coefficient vectors, the number of elements communicated per iteration increases compared to existing algorithms. Nonetheless, the total communication cost (i.e. total communicated elements) can be acceptable given i) the decreased communication round and ii) the moderate degrees of approximations in practice, see Section~\ref{subsec:DegreeApprox}. Overall, CPCA is more suitable for problems where the number of rounds outweighs the number of elements in communication.


\begin{table*}[!tb]
\centering
\caption{Comparisons of Different Distributed Nonconvex Optimization Algorithms}
\label{table:compare_result}
\begin{threeparttable}
    \begin{tabularx}{0.9\linewidth}{c *{2}{Y} c *{2}{Y}}
        \toprule
        \multirowcell{2}{{\bfseries Algorithms}} & \multicolumn{2}{c}{\bfseries Queries} & \multirowcell{2}{{\bfseries Flops}} & \multicolumn{2}{c}{{\bfseries Communication}} \\
        \cmidrule{2-3} \cmidrule{5-6}
         & $0^{\textrm{th}}$-order & $1^{\textrm{st}}$-order & & rounds & total elements \\
        \midrule
        Alg.~2 in \cite{tang2019distributed} & $\bigO{\frac{1}{\epsilon}}$\tnote{i} & / & $\bigO{\frac{F_{0}}{\epsilon}}$\tnote{ii} & $\bigO{\frac{1}{\epsilon}}$ & $\bigO{\frac{1}{\epsilon}}$ \\
        \midrule
        SONATA \cite{scutari2019distributed} & / & $\bigO{\frac{1}{\epsilon}}$ & $\bigO{\frac{\max(F_{1},F_{p})}{\epsilon}}\tnote{ii}$ & $\bigO{\frac{1}{\epsilon}}$ & $\bigO{\frac{1}{\epsilon}}$ \\
        \midrule
        {\bfseries CPCA} & $\bigO{m}$\tnote{iii} & / & $\bigOfrac{m\cdot \max(m, \log \frac{m}{\epsilon}, F_{0})}$ & $\bigO{\log \frac{m}{\epsilon}}$ & $\bigO{m \log \frac{m}{\epsilon}}$ \\
        \bottomrule
        \addlinespace[0.5ex]
    \end{tabularx}
    \begin{tablenotes}[para] 
        \item[i] $\epsilon$ is the given solution accuracy, see~\eqref{eq:sum_err_require}. 
        \item[ii] $F_{0}$, $F_{1}$, and $F_{p}$ are the costs of evaluating function values, evaluating gradients, and solving local optimization problems in \cite{scutari2019distributed}, respectively. These costs depend on the forms of objectives \cite{boyd2004convex} and hence are not explicitly specified.
        \item[iii] $m$ is the maximum degree of local approximations, and its dependence on $\epsilon$ is examined in Lemma~\ref{lem:degree_approx}, Section~\ref{subsec:DegreeApprox}.
    \end{tablenotes}
\end{threeparttable}
\end{table*}
\section{Further Discussions and Applications}\label{sec:discussion}

\subsection{Degrees of Polynomial Approximations}\label{subsec:DegreeApprox}
In the initialization stage, every agent $i$ constructs polynomial approximation $p_{i}(x)$ for $f_{i}(x)$, such that \eqref{eq:approx_require} holds. The degree $m_{i}$ of $p_{i}(x)$ depends on the specified precision $\epsilon_{1}$ and the smoothness of $f_{i}(x)$. This dependence is characterized by the following lemma.
\begin{lemma}\label{lem:degree_approx}
    If $f_{i}(x)$ and its derivatives through $f_{i}^{(v-1)}(x)$ are absolutely continuous and $f_{i}^{(v)}(x)$ is of bounded variation on $X_{i}$, then $m_{i} \sim \mathcal{O}(\epsilon_{1}^{-1/v})$. If $f_{i}(x)$ is analytic on $X_{i}$, then $m_{i} \sim \mathcal{O}(\ln\frac{1}{\epsilon_{1}})$.
\end{lemma}
\begin{proof}
    This lemma can be derived from Theorems 7.2 and 8.2 in \cite{trefethen2013approximation}. It has been proven that if $f_{i}(x)$ satisfies the differentiability condition stated in the former part of the lemma, then its degree $m$ Chebyshev interpolant $p_{i}(x)$ satisfies
    \begin{equation*}
        \max_{x \in X_{i}} |f_{i}(x) - p_{i}(x)| \leq \frac{4V}{\pi v(m-v)^{v}}.
    \end{equation*}
    If $f_{i}(x)$ analytic on $X_{i}$ is analytically continuable to the open Bernstein ellipse $E_{\rho}$, where it satisfies $|f_{i}(x)| \leq M$, then,
    \begin{equation*}
        \max_{x \in X_{i}} |f_{i}(x) - p_{i}(x)| \leq \frac{4M\rho^{-m}}{\rho - 1}.
    \end{equation*}
    By setting the above upper bounds to be less than $\epsilon_{1}$, we obtain the orders of $m_{i}$. 
\end{proof}

Lemma \ref{lem:degree_approx} implies that extremely high precision (e.g., machine epsilon) can be attained with moderate $m$ (of the order of $10^{1}\sim 10^{2}$) as long as the objective functions have certain degrees of smoothness.

\subsection{An Alternative Method of Polynomial Optimization}\label{subsec:PolyOpt}
In Section~\ref{subsec:alternative_poly_opt}, we present a strategy for optimizing the recovered approximation $p_{i}^{K}(x)$ by finding its stationary points. We now discuss an alternative method by solving semidefinite programs (SDP). This method can pave the way for further investigations of the multivariate extension, see Section~\ref{subsec:multivariate_extension} below. We offer such a reformulation, considering the coefficients of $p_{i}^{K}(x)$ given the Chebyshev polynomial basis, rather than the monomial basis\cite{blekherman2012semidefinite}. The key lies in equivalently representing a univariate polynomial by an SOS polynomial. This SOS admits a quadratic form involving a positive semidefinite matrix. This matrix appears as the optimization variable of the reformulated convex SDP.



Note that $p_{i}^{K}(x)$ is a polynomial of degree $m$ in the form of \eqref{eq:cheb_rep}, with the elements of $p_{i}^{K}=[c_{0}',\ldots,c_{m}']^{\top}$ being its Chebyshev coefficients. To simplify the notation, we utilize the translation and scaling of $p_{i}^{K}(x)$ defined on $[a,b]$ to obtain
\begin{equation*}
    g_{i}^{K}(u) \triangleq p_{i}^{K}\Big(\frac{b-a}{2}u + \frac{a+b}{2}\Big)  = \sum_{j=0}^{m} c_{j}'T_{j}(u), ~~ u \in [-1,1].
\end{equation*}
Then, the optimal value of $p_{i}^{K}(x)$ on $[a,b]$ and that of $g_{i}^{K}(u)$ on $[-1,1]$ are equal, and the optimal points $x_{p}^{*}$ and $u_{g}^{*}$ satisfy
\begin{equation}\label{eq:opt_point_relation}
    x_{p}^{*} = \frac{b-a}{2} u_{g}^{*} + \frac{a+b}{2}.
\end{equation}
Therefore, once we have solved the following problem
\begin{equation}\label{problem:poly_main}
    \min_{u}~g_{i}^{K}(u) \quad \textrm{s.t. } u\in [-1,1],
\end{equation}
we can use \eqref{eq:opt_point_relation} to obtain the optimal value and optimal points of $p_{i}^{K}(x)$ on $[a,b]$. We now discuss how to transform problem \eqref{problem:poly_main} to a convex optimization problem.

We first transform problem \eqref{problem:poly_main} to its equivalent form
\begin{equation}\label{problem:poly_epigraph}
    \max_{t}~t \qquad \textrm{s.t. } g_{i}^{K}(x)-t \geq 0, ~~ \forall x \in [-1,1].
\end{equation}
The equivalence follows from the fact that $(x^{*},t^{*})$ is optimal for problem \eqref{problem:poly_epigraph} \textit{iff} $x^{*}$ is optimal for problem \eqref{problem:poly_main} and $t^{*} = g_{i}^{K}(x^{*})$. To further transform the inequality constraint, we utilize the non-negativity of $g_{i}^{K}(x)-t$ for $x \in [-1,1]$, where $t$ is viewed as a specified parameter. It follows from the \textit{Markov-Luk{\'a}cs theorem}\cite[Theorem 3.72]{blekherman2012semidefinite} that $g_{i}^{K}(x)-t$ is non-negative for $x\in [-1,1]$ \textit{iff} it can be expressed as
\begin{equation*}  
    g_{i}^{K}(x)-t =\left\{
    \begin{aligned}
        &(x+1)h_{1}^{2}(x) + (1-x)h_{2}^{2}(x), ~~ \textrm{if $m$ is odd}, \\
        &h_{1}^{2}(x) + (x+1)(1-x)h_{2}^{2}(x), ~~ \textrm{if $m$ is even},
    \end{aligned}
    \right.
\end{equation*}
where $h_{1}(x)$ and $h_{2}(x)$ are polynomials of degree $\lfloor \frac{m}{2} \rfloor$ and $\lfloor \frac{m-1}{2} \rfloor$, respectively, and some coefficients of $h_{1}(x)$ and $h_{2}(x)$ depend on $t$. Based on \eqref{eq:sos_psd}, there exist $Q,Q'\in \mathbb{S}_{+}$ such that
\begin{equation*}
    h_{1}^{2}(x) = v_{1}(x)^{\top}Qv_{1}(x), \quad h_{2}^{2}(x) = v_{2}(x)^{\top}Q'v_{2}(x),
\end{equation*}
where
\begin{alignat*}{2}
    v_{1}(x) &= [T_{0}(x),\ldots,T_{d_{1}}(x)]^{\top}, &\quad& d_{1} = {\textstyle \lfloor \frac{m}{2} \rfloor}, \\
    v_{2}(x) &= [T_{0}(x),\ldots,T_{d_{2}}(x)]^{\top}, && d_{2} = {\textstyle \lfloor \frac{m-1}{2} \rfloor}.
\end{alignat*}
Based on the above transformations, by ensuring that the Chebyshev coefficients of $g_{i}^{K}(x)-t$ are consistent, we transform the inequality constraint in problem \eqref{problem:poly_epigraph} to a set of equality constraints related to the elements of $Q$ and $Q'$. When $m$ is odd, we reformulate problem \eqref{problem:poly_epigraph} as
\begin{align}\label{problem:odd_reformulated}
    \max_{t,Q,Q'} \quad & t \notag \\
    \textrm{s.t.} \quad & c_{0}' = t + Q_{00} + Q'_{00} + \frac{1}{2} \Big(\sum_{u=1}^{d_{1}+1}Q_{uu} + \sum_{u=1}^{d_{2}+1} Q'_{uu} \Big) \notag \\
              &\qquad + \frac{1}{4} \sum_{|u-v|=1} \left(Q_{uv} - Q'_{uv}\right), \notag  \\
        & c_{j}' = \frac{1}{2} \sum_{(u,v)\in \mathcal{A}} \left(Q_{uv} \!+\! Q'_{uv} \right) \!+\! \frac{1}{4} \sum_{(u,v)\in \mathcal{B}} \left(Q_{uv} \!-\! Q'_{uv}\right), \notag \\
        & \qquad j=1,\ldots,m, \notag \\
        & Q \in \mathbb{S}_{+}^{d_{1} +1}, ~ Q' \in \mathbb{S}_{+}^{d_{2} +1},
\end{align}
where the rows and columns of $Q \textrm{ and } Q'$ are indexed by $0,1,\ldots$, and
\begin{subequations}\label{eq:set_A_B}
\begin{align}
    \mathcal{A} &= \{(u,v)|u+v=i \vee |u-v|=i\}, \label{eq:set_A}  \\
    \mathcal{B} &= \big\{(u,v)|u+v=i-1 \vee |u-v|=i-1 \vee \nonumber \\
                &\phantom{=\{(u,v)|~} |u+v-1|=i \vee \big||u-v|-1\big|=i \big\}. \label{eq:set_B}
\end{align}
\end{subequations}
When $m$ is even, we reformulate problem \eqref{problem:poly_epigraph} as
\begin{align}\label{problem:even_reformulated}
    \max_{t,Q,Q'} \quad & t \notag \\
    \textrm{s.t.} \quad & c_{0}' = t + Q_{00} + \frac{1}{2} Q'_{00} + \frac{1}{2}\sum_{u=1}^{d_{1}+1}Q_{uu} \notag \\
        &\qquad + \frac{1}{4} \sum_{u=1}^{d_{2}+1} Q'_{uu} + \frac{1}{8} \sum_{|u-v|=2} Q'_{uv}, \notag \\
        & c_{j}' = \frac{1}{2} \sum_{(u,v)\in \mathcal{A}} \Big(Q_{uv} + \frac{1}{2} Q'_{uv} \Big) + \frac{1}{8} \sum_{(u,v)\in \mathcal{C}} Q'_{uv}, \notag \\
        & \qquad j=1,\ldots,m, \notag \\ 
        & Q \in \mathbb{S}_{+}^{d_{1} +1},~Q' \in \mathbb{S}_{+}^{d_{2} +1},
\end{align}
where $\mathcal{A}$ is given by \eqref{eq:set_A}, and
\begin{equation*}
    \begin{split}
        \mathcal{C} = \big\{(u,v)|&u+v=i-2 \vee |u-v|=i-2 \\
                    & \vee |u+v-2|=i \vee \big||u-v|-2\big|=i \big\}.
    \end{split}
\end{equation*}
The following theorem guarantees the equivalence of these reformulations and problem \eqref{problem:poly_main}.
\begin{theorem}\label{thm:sdp_equivalence}
    When $m$ is odd (resp., even), problem \eqref{problem:odd_reformulated} (resp., problem \eqref{problem:even_reformulated}) has the same optimal value as problem \eqref{problem:poly_main}.
\end{theorem}
\begin{proof}
    See Appendix \ref{subsec:appendix_derivation}. 
\end{proof}

Both of these reformulated problems are SDPs and, therefore, can be efficiently solved via primal-dual interior-point methods\cite{boyd2004convex}. The solving process is terminated when
\begin{equation*}
    \mathbf{G_3'}: \quad 0\leq f_{e}^{*} - p^{*} \leq \epsilon_{3}
\end{equation*}
is achieved, where $f_{e}^{*}$ is the returned estimate of the optimal value $p^{*}$ of $p_{i}^{K}(x)$ on $X=[a,b]$, and $\epsilon_{3}>0$ is some specified tolerance. In this case, the requirement \eqref{eq:sum_err_require} becomes $\sum_{i=1}^{3} \epsilon_i = \epsilon, \epsilon_i > 0$. The optimal points of $g_{i}^{K}(x)$ can be obtained by i) plugging in $f_e^*$ and solving the polynomial root-finding problem $g_i^K(x) = f_e^*$, or ii) using the complementary slackness condition to obtain a measure with support in the optimal points\cite{blekherman2012semidefinite}. We can then calculate the optimal points of $p_{i}^{K}(x)$ on $X$ by \eqref{eq:opt_point_relation}.

\subsection{Application Scenarios}

\emph{Distributed optimization with decoupling}:
In some cases, distributed optimization problems with high-dimensional variables are naturally decoupled or can be converted to decoupled via a change of coordinates or variables. Then, we can separately optimize over the single variable to solve the whole problem. 
For example, consider a regularized facility location problem\cite{boyd2004convex,metel2019simple}. Every agent $i(=1,\ldots,N)$ owns a location $u_{i} \in \mathbb{R}^{n}$ and aims to agree on a point $x \in \mathbb{R}^{n}$ that minimizes the sum of the weighted taxicab distances to all the locations and a regularizer, i.e.,
\begin{equation}\label{problem:facility_location} 
    \min_{x \in X^n} ~ \sum_{i=1}^{N} w_{i} \|x - u_{i}\|_{1} + g(x),
\end{equation}
where $X^n \subset \mathbb{R}^n$ is a constraint set, $w_{i} \geq 0$ is a weight, $g(x) \triangleq \sum_{i=1}^{n} \kappa \log(1+|x(k)|/\nu)$ is the log-sum penalty, and $\kappa,\nu>0$ inside $g(x)$ are parameters. Note that $g(x)$ is a common surrogate of the $\ell_0$-regularizer to promote sparsity, and it is nonconvex and Lipschitz continuous\cite{metel2019simple}. We can solve problem~\eqref{problem:facility_location} by considering the following set of subproblems with univariate nonconvex Lipschitz objectives for $k = 1,\ldots,n$
\begin{equation*}
    \min_{x(k) \in X} ~ \sum_{i=1}^{N} w_{i} |x(k) - u_{i}(k)| + \kappa \log(1+|x(k)|/\nu).
\end{equation*}



\emph{Hyperparameter optimization in distributed learning}:
We want to find a continuous hyperparameter (e.g., learning rate or regularization parameter) that leads to an optimal model with the minimum generalization error\cite{yang2020hyperparameter,gu2021optimizing}. This problem is 
\begin{equation*}
 \begin{split}
     \min_{\lambda} ~& F(\lambda) = \frac{1}{N} \sum_{i=1}^{N} R_i(w(\lambda),\lambda) \\
         \text{s.t.} ~& w(\lambda) = \argmin_{w} \frac{1}{N} \sum_{i=1}^{N} L_i(w,\lambda),
 \end{split}       
\end{equation*}
where $\lambda$ is the hyperparameter, $w(\lambda)$ are the model parameters, $L_i(w,\lambda)$ is the loss function of agent $i$ in the learning problem (e.g., regularized empirical risk on the local training dataset), and $R_i(w(\lambda),\lambda)$ is the generalization error of the learned model given $\lambda$ (e.g., empirical risk on the local validation dataset). In practice, we usually obtain approximate model parameters $\tilde{w}(\lambda)$ via some descent-based algorithms $\mathcal{A}$, i.e.,
    $w(\lambda) \approx \tilde{w}(\lambda) = \mathcal{A}(\lambda)$.
Based on \cite[Theorem~1]{gu2021optimizing}, if $\lambda$ is a continuous variable and $\mathcal{A}(\lambda)$ and $R(w(\lambda),\lambda)$ are continuous functions, then $F(\lambda)$ is a univariate continuous function of $\lambda$. In this case, the coupling in $w(\lambda)$ may pose additional challenges. The proposed algorithm needs further adjustment to fit this interesting problem, which is of independent interest in the learning community. 




\textit{Distributed data statistics}:
In a networked system, agents hold a portion of the whole data. The goal is to estimate global data statistics in a distributed manner. Typical examples include \romannumeral1) estimating a global probability density function (PDF), which is a univariate function for a specific random variable, and \romannumeral2) estimating a separable functional $T(P)$ (e.g., entropy) of an unknown discrete distribution $P=(p_1,p_2,\ldots)$\cite{wu2021polynomial}, where $T(P) = \sum_{j} f(p_j)$
for some univariate function $f$. Instead of directly transmitting plenty of local data, agents can solve problems with a low communication cost by exploiting approximations of local PDFs or functions $f$ and utilizing consensus-based fusion of local information. 

\subsection{The Multivariate Extension}\label{subsec:multivariate_extension}
We discuss the design and open issues of the multivariate extension of the proposed algorithm. The differences lie in the stages of initialization and local optimization of approximations. Specifically, any smooth local objective $f_{i}(x)$ in the $n$-dimensional hypercube $[-1,1]^{n}$ can be approximated by a multivariate Chebyshev polynomial of the following form 
\begin{equation*}
    \hat{f}_i(x) = \sum_{k_1=0}^{m_1}\cdots\sum_{k_n=0}^{m_n} a_{k_1,\ldots,k_s} T_{k_1}(x(1))\cdots T_{k_n}(x(n)),
\end{equation*}
where $x(n)$ denotes the $n$-th element in a vector $x$, and $T_{k_n}(\cdot)$ denotes the $k_n$-th Chebyshev polynomial defined on $[-1,1]$\cite{trefethen2017multivariate}. After local approximations $\hat{f}_{i}(x)$ are constructed, agents exchange and update their local variables that store those coefficients (as discussed in Section~\ref{subsec:iterations}) and obtain an approximation of the global objective function. Finally, they locally optimize this approximation by solving a hierarchy of SDPs (of increasing size, which is related to degrees of polynomials used in the SOS representation) to obtain a monotone non-decreasing sequence of lower bounds that asymptotically approach the globally optimal value.

There remain two technical issues. The first is to find a systematic means of constructing $\hat{f}_i(x)$ of an appropriate degree that strictly meets the specified accuracy requirement \eqref{eq:approx_require}. The second is to reduce the costs of flops in local optimization due to hierarchical SDPs. 
Nonetheless, the central idea of introducing polynomial approximation can still \romannumeral1) offer a new perspective of representing complicated objective functions by simple vectors of coefficients of approximations, which allows the dissemination of the global objective in networked systems, and \romannumeral2) facilitate obtaining approximate globally optimal solutions in a distributed and asymptotic manner. Thanks to the modular structure of the algorithm, the new advances in polynomial approximation can be incorporated to better find approximate globally optimal solutions.
\section{Numerical Evaluations}\label{sec:experiment}
We present simulation results to illustrate the performance of CPCA and compare it with other algorithms.
We generate a network with $N=30$ agents using the Erd\H{o}s-R{\'e}nyi model with connectivity probability $0.4$. In other words, the probability that every pair of agents can communicate is $0.4$. Suppose that all the local constraint sets are the same interval $X=[-1,1]$. We consider two instances of problem \eqref{problem:main_focus} with different types of local objective functions. In the first instance, the local objective of agent $i$ is
\begin{equation}\label{eq:local_obj_exp}
    f_{i}(x) = a_{i}e^{b_{i}x} + c_{i}e^{-d_{i}x},
\end{equation}
where $a_{i},b_{i}\sim U(1,2),c_{i},d_{i}\sim U(2,4)$ are uniformly distributed. Hence, $f_{i}(x)$ is convex and Lipschitz continuous on $X$. In the second instance, the local objective of agent $i$ is 
\begin{equation}\label{eq:local_obj_log_exp}
    f_{i}(x) = \frac{a_{i}}{1+e^{-x}} + b_{i}\log (1+x^{2}),
\end{equation}
where $a_{i}\sim \mathcal{N}(10,2),b_{i}\sim \mathcal{N}(5,1)$ are Gaussian random variables. 
In this case, $f_{i}(x)$ is nonconvex and Lipschitz continuous on $X$. We use the Chebfun toolbox\cite{driscoll2014chebfun} that implements the chopping technique\cite{aurentz2017chopping} to help construct Chebyshev polynomial approximations.

For comparison, we implement SONATA-L\cite{scutari2019distributed}, Alg.~2 in \cite{tang2019distributed}, projected distributed sub-gradient descent method (Proj-DGD, for convex objectives)\cite{nedic2010constrained}, and projected stochastic gradient descent method (Proj-SGD, for nonconvex objectives)\cite{bianchi2012convergence}. The step sizes of SONATA-L, Alg.~2 in \cite{tang2019distributed}, Proj-DGD, and Proj-SGD are set as $\alpha^{t} = \alpha^{t-1}(1- 0.01 \alpha^{t-1})$ with $\alpha^{0}=0.5$, $\alpha^{t} = 0.05$, $\alpha^{t} = 1/t^{0.5}$, and $\alpha^{t}=0.1/t^{0.9}$, respectively, based on the guidelines therein. For Alg.~2 in \cite{tang2019distributed}, the smoothing radius $u^{t}$ present in gradient estimates is given by $u^{t}=1/t^{3/4}$. For Proj-SGD, we consider $50$ Monte-Carlo runs and plot the curves of average objective errors.

\begin{figure}[t]
\vspace{-1ex}
\centering
  \subfloat[Communication]{\includegraphics[width=0.5\columnwidth]{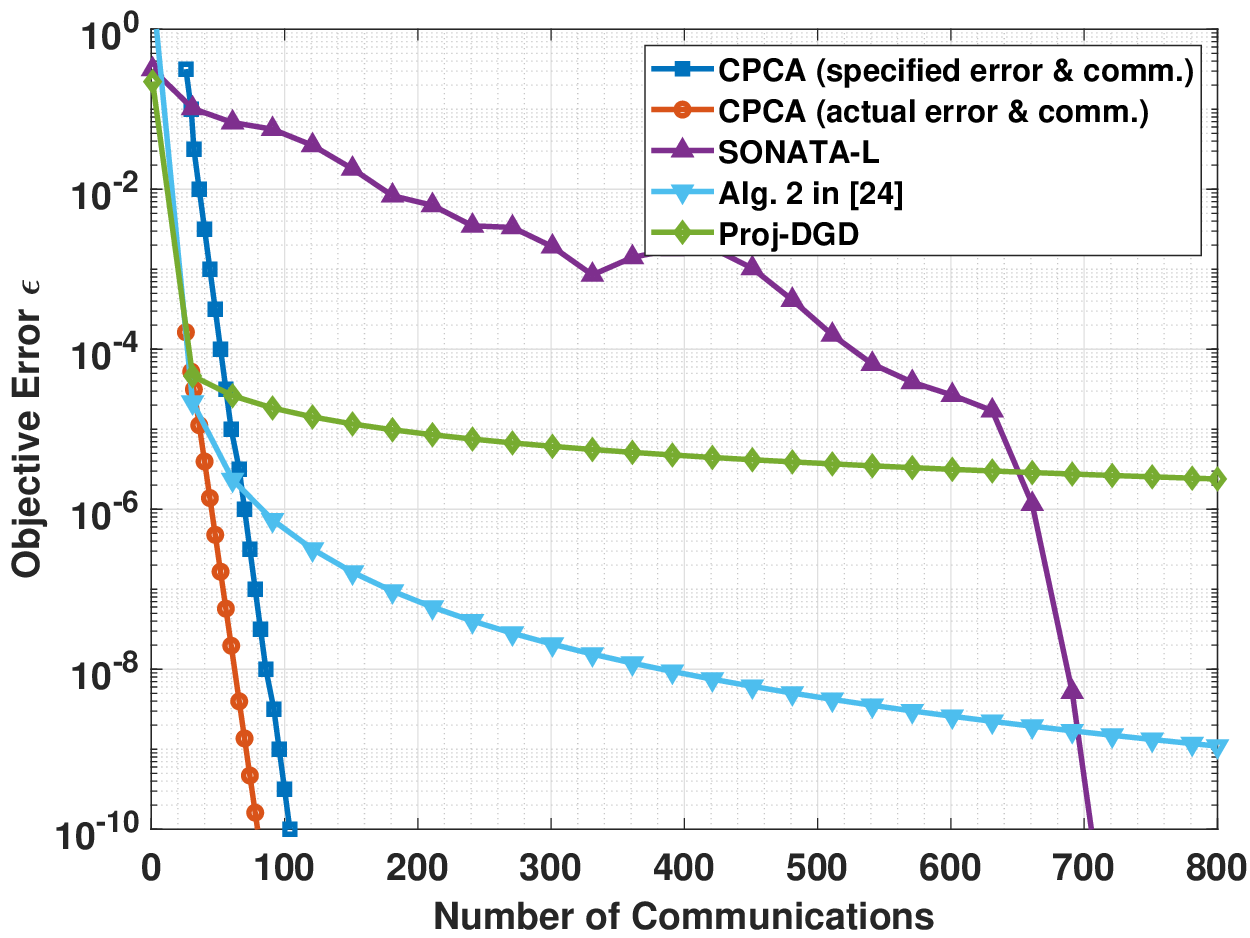}\label{fig:cvx_cmmunctn}} \hfil
  \subfloat[Oracle Queries]{\includegraphics[width=0.5\columnwidth]{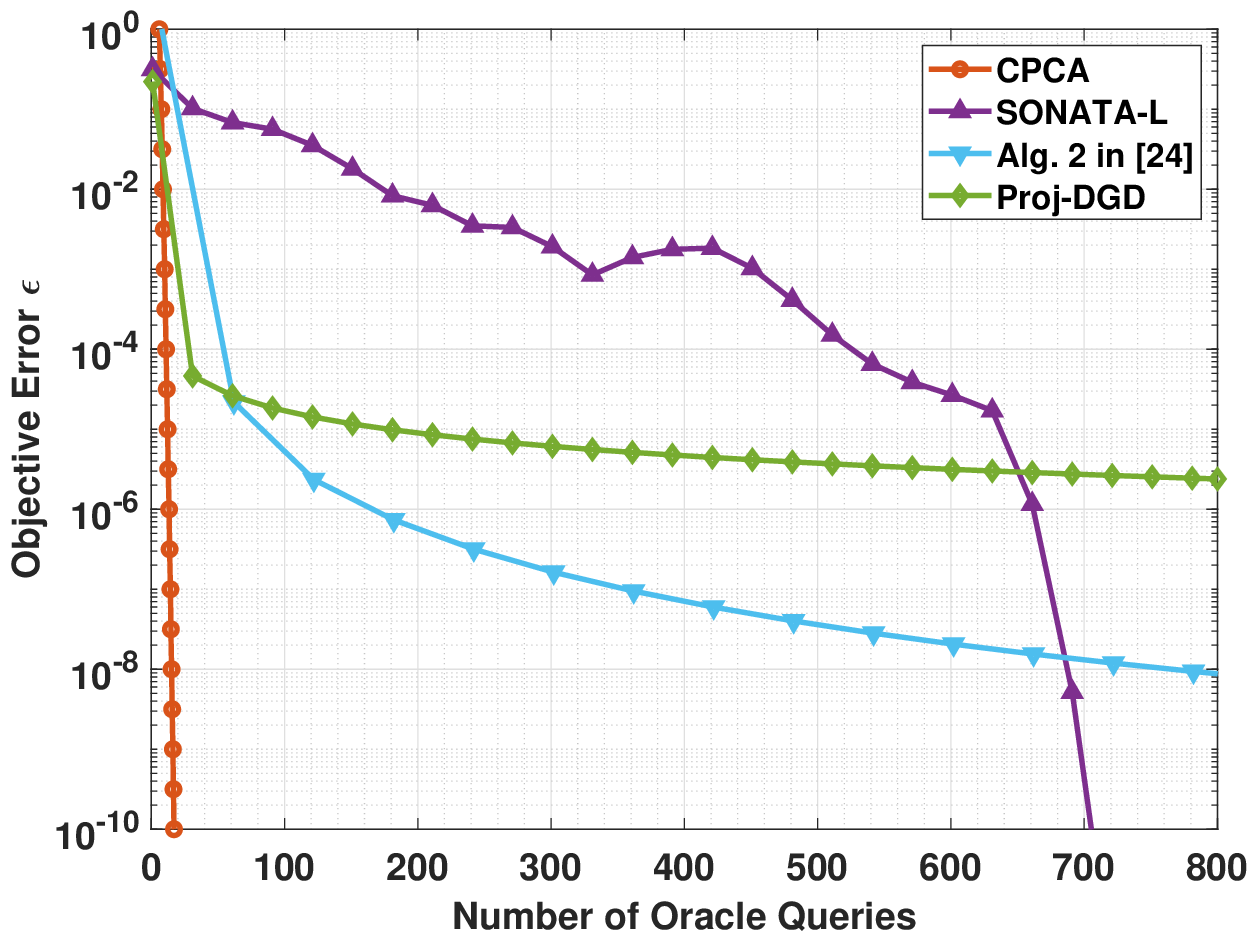}\label{fig:cvx_orclqur}}
  \caption{Comparison of CPCA, SONATA-L, Alg.~2 in \cite{tang2019distributed}, and Proj-DGD for solving problem \eqref{problem:main_focus} with \eqref{eq:local_obj_exp} as local objectives.} 
  \label{fig:cvxComplexity} 
  \vspace{-1.5em}
\end{figure}

\begin{figure}[t]
\centering
  \subfloat[Communication]{\includegraphics[width=0.5\columnwidth]{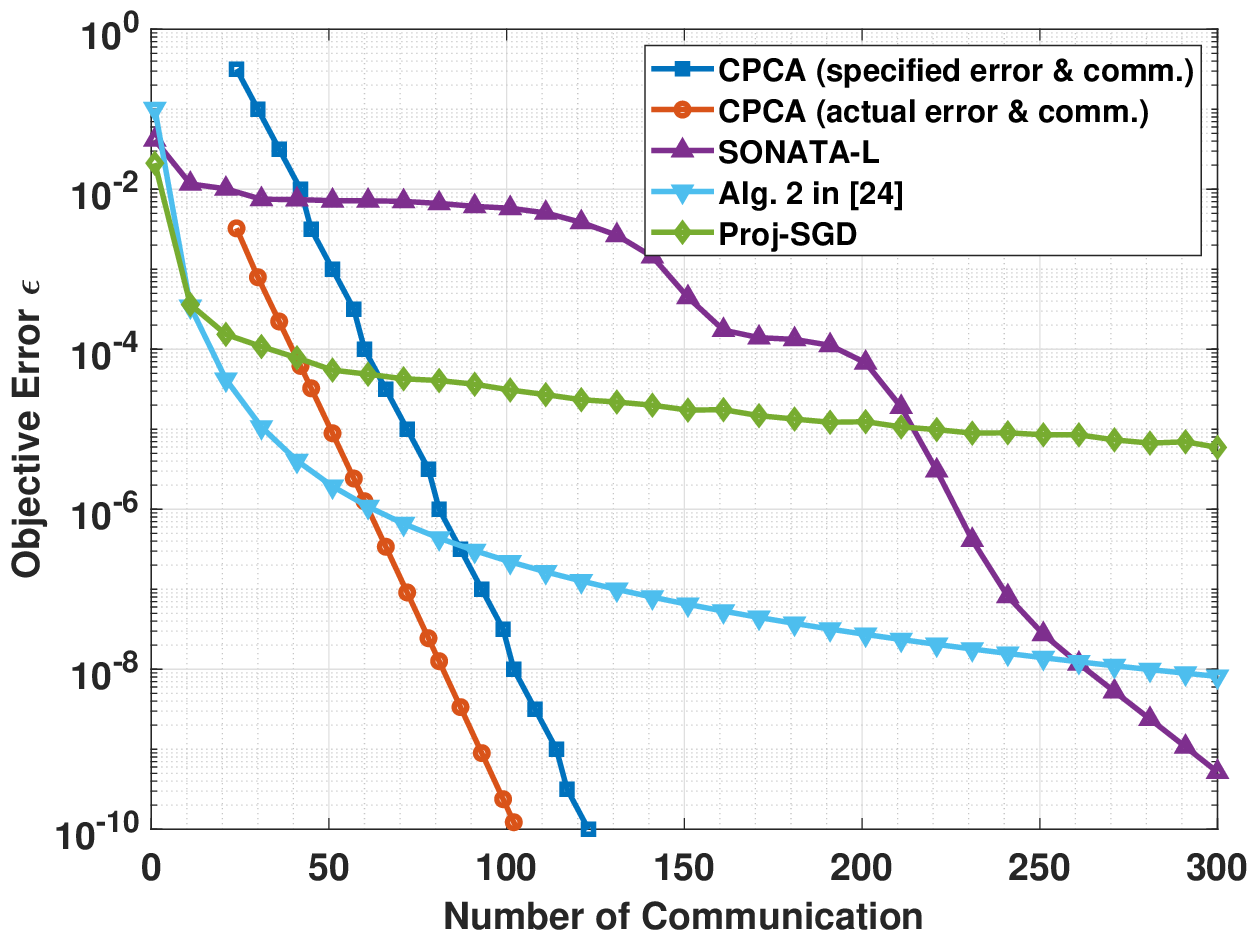}\label{fig:cmmunctn}} \hfil
  \subfloat[Oracle Queries]{\includegraphics[width=0.5\columnwidth]{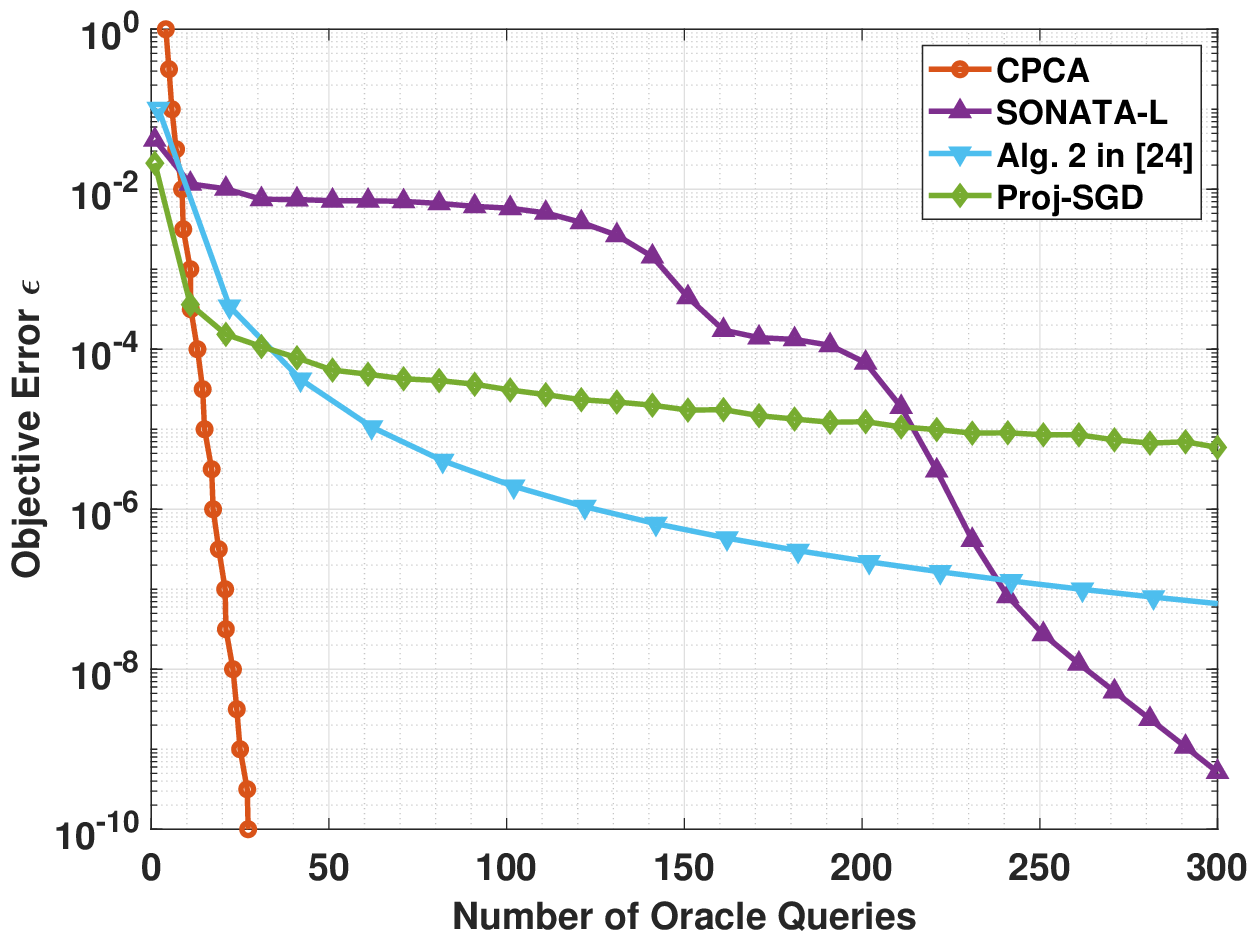}\label{fig:orclqur}}
  \caption{Comparison of CPCA, SONATA-L, Alg.~2 in \cite{tang2019distributed}, and Proj-SGD for solving problem \eqref{problem:main_focus} with \eqref{eq:local_obj_log_exp} as local objectives.} 
  \label{fig:ncvxComplexity}
\end{figure}

Fig.~\subref*{fig:cvx_cmmunctn} and \subref*{fig:cmmunctn} show the relationships between the objective error $\epsilon$ and the number of inter-agent communication rounds. For CPCA and other algorithms, $\epsilon$ denote $|f_{e}^{*}-f^{*}|$ and $|f(\overbar{x}^{t})-f^{*}|$, respectively, where $\overbar{x}^{t}$ is the average of all agents' local estimates at time $t$. Specifically, the blue line shows how the specified accuracy influences the number of communication rounds executed by CPCA, whereas the orange line indicates the actual objective error after the algorithm performs communication and obtains solutions.
The orange line is below the blue line (i.e., the actual error is less than the specified error). It stems from the careful design of CPCA to satisfy the given accuracy.
Moreover, to reach a certain precision, CPCA requires a less round of communication thanks to its linearly convergent inner iterations. However, it communicates more elements per iteration. For example, to reach a precision of $10^{-10}$, CPCA exchanges vectors of size $19$ and $29$ when the local objective function is \eqref{eq:local_obj_exp} and \eqref{eq:local_obj_log_exp}, respectively. In comparison, the remaining algorithms exchange scalars of size one.

Fig.~\subref*{fig:cvx_orclqur} and \subref*{fig:orclqur} show the relationships between the objective error $\epsilon$ and the number of oracle queries. For CPCA that requires zeroth-order queries, the horizontal axis represents the average number of queries for one agent (i.e., the average degree of local proxies plus one). Note that such queries are not performed during consensus-based iterations. These results support the discussions in Section~\ref{subsec:DegreeApprox} that $m$ depends on $\epsilon$, and that extremely small $\epsilon$ is generally associated with moderate $m$. 
Furthermore, because the size of the coefficient vector exchanged per iteration equals the number of sampling points (i.e., the required zeroth-order queries), these results also indicate the per-step communication costs of CPCA\@.
Alg.~2 in \cite{tang2019distributed} and the remaining algorithms require two zeroth-order queries and one first-order query at every iteration, respectively. Hence, the curves corresponding to other algorithms except \cite{tang2019distributed} are identical to those in Fig.~\subref*{fig:cvx_cmmunctn} and \subref*{fig:cmmunctn}. We observe that CPCA requires fewer oracle queries in both examples. This difference results from its design of gradient-free iterations. 

\begin{figure}[t]
    \vspace{-1ex}
    \centering
    \subfloat[Cycle graph]{\includegraphics[width=0.5\linewidth]{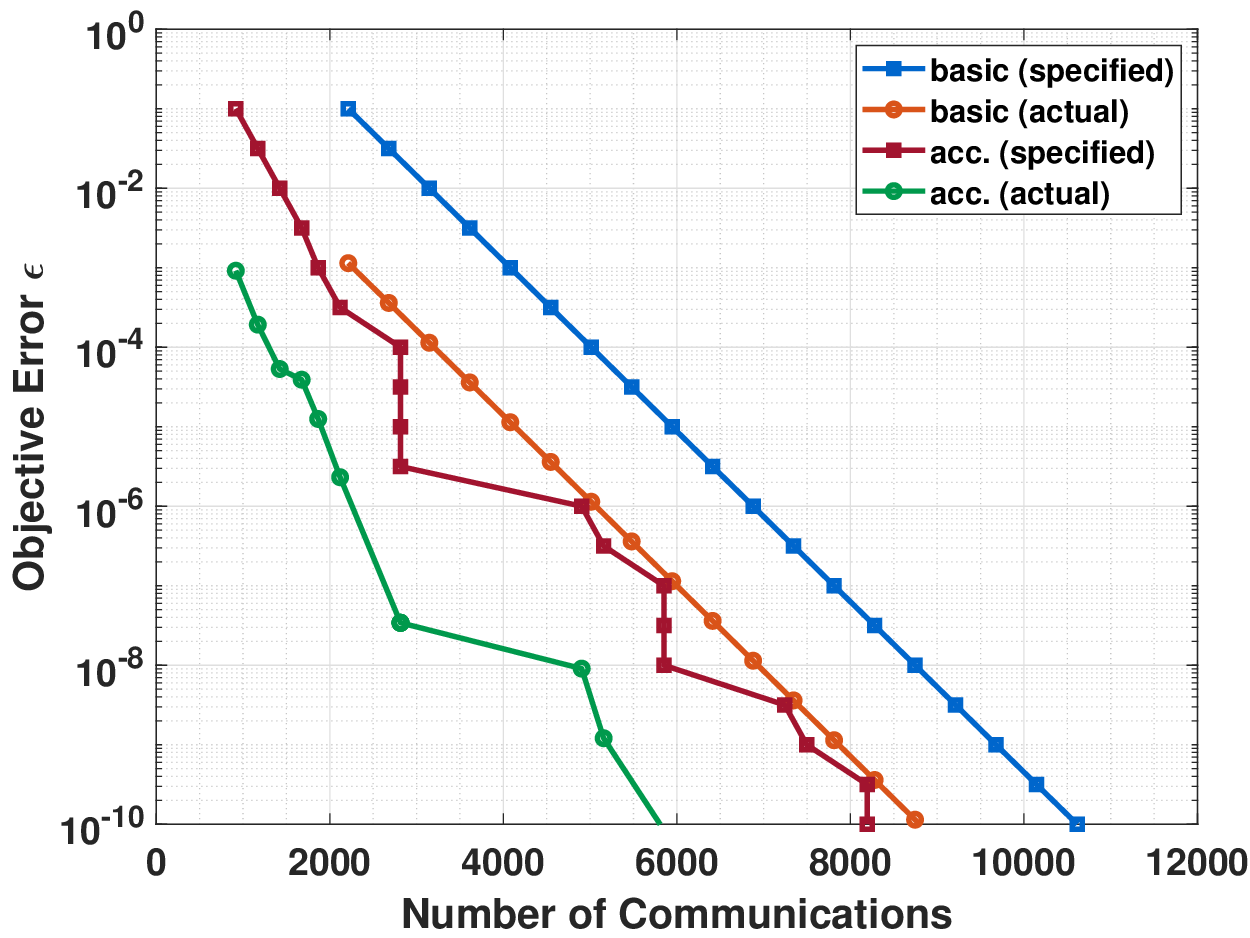}\label{fig:cpca_cycle}} \hfil
    \subfloat[Random graph]{\includegraphics[width=0.5\linewidth]{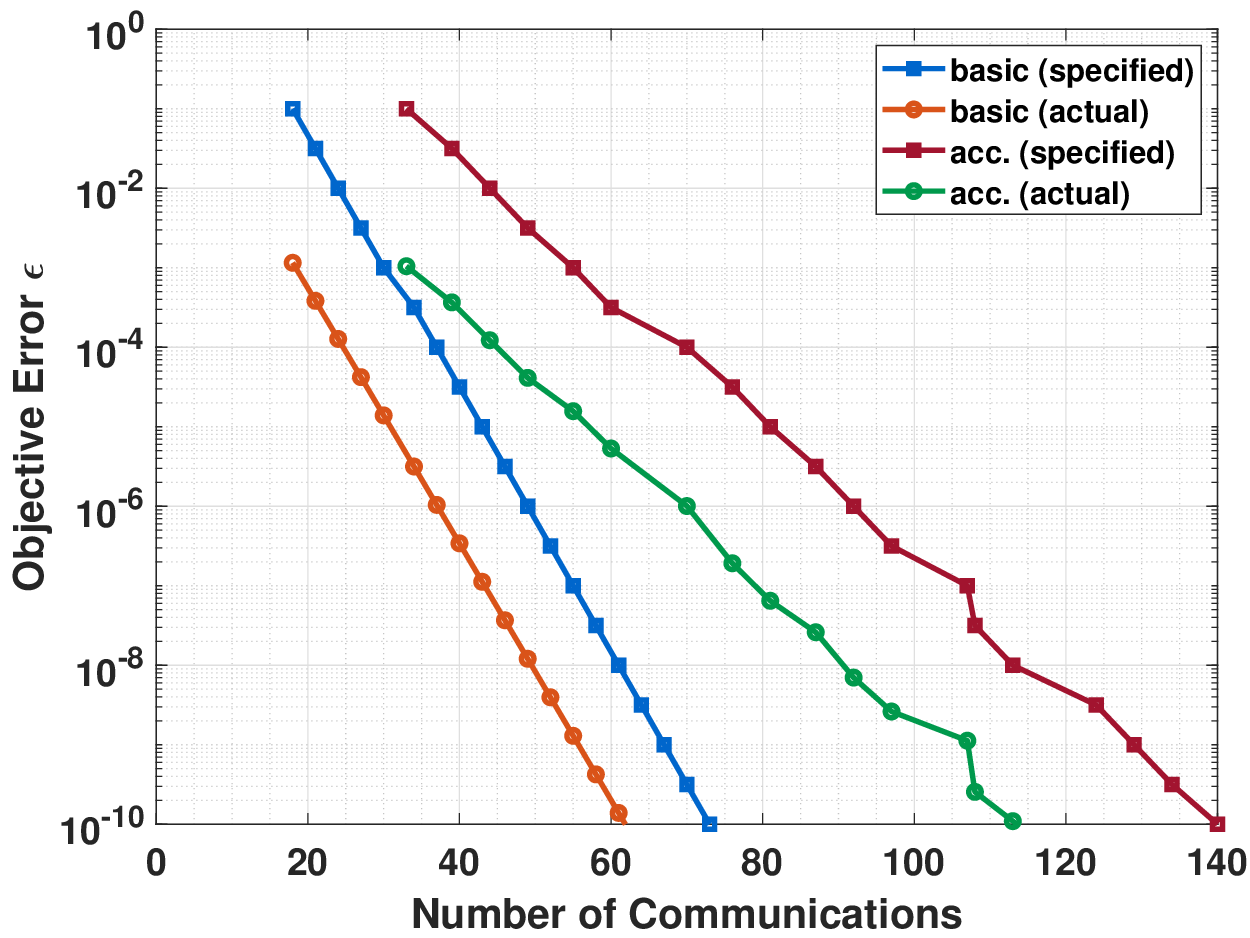}\label{fig:cpca_random}}
    \caption{Performance of the proposed algorithm with the basic and accelerated (corresponding to ``acc.'' in the captions) consensus methods.}
    \label{fig:CPCA_performance} 
\end{figure}

We further analyze how the network topology and consensus-based iterations influence the performance. We generate two networks with $N=100$ agents. The first one is a cycle graph, where every agent can only communicate with its two direct neighbors. The second one is an Erd\H{o}s-R{\'e}nyi graph with connectivity probability $0.4$. We implement CPCA with the basic \eqref{eq:avg_consensus_implement} and accelerated consensus methods \cite{olshevsky2017linear} (adopted in the conference version\cite{he2020cpca}) to solve problem~\eqref{eq:local_obj_log_exp}. Here we allow centralized stopping of consensus iterations. The relationships between the objective error and the number of inter-agent communication rounds are shown in Fig.~\ref{fig:CPCA_performance}. 
The blue and red lines show how the specified accuracy determines the number of communication rounds executed by CPCA. The orange and green lines indicate the relationship between the number of communication rounds and the attained objective error. We observe that the algorithmic performance with different consensus methods depends on the network topology. On this cycle graph, the accelerated method wins, whereas on this random graph, the basic method wins.

\section{Conclusion}\label{sec:conclusion}
We proposed CPCA to solve distributed optimization problems with nonconvex Lipschitz continuous univariate objectives and convex local constraint sets. The proposed algorithm leverages Chebyshev polynomial approximation, consensus, and polynomial optimization. It consists of three stages, i.e., \romannumeral1) constructing approximations of local objectives, \romannumeral2) performing average consensus to disseminate information, and \romannumeral3) locally optimizing the obtained approximation of the global objective via polynomial optimization techniques. We provided comprehensive theoretical analysis and numerical results to illustrate its effectiveness, including \romannumeral1) obtaining $\epsilon$-globally optimal solutions, \romannumeral2) being efficient in oracle and communication complexities, and \romannumeral3) achieving distributed termination. Future directions include but are not limited to designing efficient multivariate extensions via polynomial approximations and handling noisy evaluations of function values. 

\appendix
\section{Appendix}
\subsection{Proof of Theorem~\ref{thm:itr_precision}}\label{subsec:proof_itr}
    We rewrite the update \eqref{eq:avg_consensus_implement} as
    $p_{i}^{t+1} = \sum_{j=1}^{N} w_{ij} p_{j}^{t}$, where
    \begin{equation*}
        w_{ij} = \left\{
        \begin{aligned}
            & (2\mathrm{max}(d_{i},d_{j}))^{-1}, \quad & j\in \mathcal{N}_{i}, j\neq i, \\
            & 0, & j \notin \mathcal{N}_{i}, j\neq i, \\
            & 1-{\textstyle \sum_{j \in \mathcal{N}_{i}}} w_{ij}, & j = i.
        \end{aligned}
        \right.
    \end{equation*}
    It follows that $W \triangleq \left(w_{ij}\right)_{i,j=1}^{N}$ is row stochastic, i.e., $\sum_{j=1}^{N} w_{ij} = 1, 0 \leq w_{ij} \leq 1, \forall i,j = 1,\ldots,N$.
    We consider the $k$-th element of $p_i^t$, where $k=0,\ldots,m$ and $i\in \mathcal{V}$. Using the row-stochasticity of $W$, we have
    \begin{equation*}
        \begin{split}
            p_{i}^{t+1}(k) = \sum_{j=1}^{N} w_{ij} p_{j}^{t}(k) \leq \sum_{j=1}^{N} w_{ij} \max_{j\in \mathcal{V}} p_{j}^{t}(k) = \max_{j\in \mathcal{V}} p_{j}^{t}(k). 
        \end{split}
    \end{equation*}
    Let $\max_{j\in \mathcal{V}} p_{j}^{t}(k) \triangleq M^{t}(k)$ and $\min_{j\in \mathcal{V}} p_{j}^{t}(k) \triangleq m^{t}(k)$. Hence, $M^{t+1}(k)\leq M^{t}(k), m^{t+1}(k)\geq m^{t}(k)$.
    The convergence of \eqref{eq:avg_consensus} implies that $\lim_{t\rightarrow \infty} p_{i}^{t}(k) = \overbar{p}(k), \forall i\in \mathcal{V}$. Therefore, $\lim_{t\rightarrow \infty} M^{t}(k) = \overbar{p}(k)$. Given that $M^{t}(k)$ is non-increasing with respect to $t$, we have $\overbar{p}(k)\leq M^{t}(k), \forall t\in \mathbb{N}$. We apply a similar reasoning for $m^t(k)$ and obtain
        $m^{t}(k) \leq \overbar{p}(k) \leq M^{t}(k), \forall t\in \mathbb{N}$.
    Suppose that agents terminate at the $K$-th iteration. Since the max/min consensus is guaranteed to converge within $U$ iterations, we have
    \begin{equation*}
        r_{i}^{K}(k)-s_{i}^{K}(k) = M^{K'}(k)-m^{K'}(k),
    \end{equation*}
    where $K'\triangleq K-U$. The stopping criterion $\|r_{i}^{K} - s_{i}^{K}\|_{\infty} \leq \delta$ is equivalent to $r_{i}^{K}(k)-s_{i}^{K}(k)\leq \delta, \forall k$. When it is satisfied, 
    \begin{equation*}
        \left|p_{i}^{K}(k)\!-\!\overbar{p}(k)\right| \!\leq\! M^{K}(k)\!-\!m^{K}(k) \!\leq\! r_{i}^{K}(k)\!-\!s_{i}^{K}(k) \!\leq\! \delta, ~\forall i,k.
    \end{equation*}

\subsection{Proof of Theorem~\ref{thm:alg_accuracy_result}}\label{subsec:proof_accuracy}

    First, we establish the closeness between $p_{i}^{K}(x)$ and $\overbar{p}(x)$. Theorem~\ref{thm:itr_precision} implies that $\|p_{i}^{K} - \overbar{p}\|_{\infty} \leq \delta, \forall i\in \mathcal{V}$.
    Suppose that the sets of Chebyshev coefficients of $p_{i}^{K}(x)$ and $\overbar{p}(x)$ are $\{c_{j}'\}$ and $\{\overbar{c}_{j}'\}$, respectively. It follows that
        $\left|c_{j}'-\overbar{c}_{j}'\right| \leq \delta, \forall j=0,\ldots,m$.
    Consequently,
    \begin{align*}
        \begin{split}
        |p_{i}^{K}(x) - \overbar{p}(x)| &= \bigg| \sum_{j=0}^{m} (c_{j}'-\overbar{c}_{j}') T_{j} \left(\frac{2x-(a+b)}{b-a}\right) \bigg| \\
        &\leq \sum_{j=0}^{m} \left|c_{j}' - \overbar{c}_{j}'\right| \cdot 1 \leq \delta(m+1) = \epsilon_{2},
        \end{split}
    \end{align*}
    where we use the fact that $\left|T_{j}(x)\right|\leq 1, \forall x \in [-1,1]$.


    Then, we establish the closeness between $\overbar{p}(x)$ and $f(x)$. Since $\overbar{p} $ is the average of all $p_{i}^{0}$, $\overbar{p}(x)$ is also the average of all $p_{i}(x)$. Based on the results in Section~\ref{subsec:init}, we have
    \begin{equation*}
        \begin{split}
        |\overbar{p}(x) - f(x)| &= \bigg|\frac{1}{N}\sum_{i=1}^{N} \big(p_{i}(x) - f_i(x)\big)\bigg| \\
         &\leq \frac{1}{N} \sum_{i=1}^{N} \left|p_{i}(x) - f_i(x)\right| \leq \frac{1}{N} N \epsilon_{1} = \epsilon_{1}.
        \end{split}
    \end{equation*}
    Note that $\epsilon_{1} = \epsilon_{2} = \epsilon/2$. Hence,
    \begin{align*}
        \left|p_{i}^{K}(x) \!-\! f(x)\right| & \leq \left|p_{i}^{K}(x) \!-\! \overbar{p}(x)\right| \!+\! \left|\overbar{p}(x) \!-\! f(x)\right|  \leq \epsilon_{1} \!+\! \epsilon_{2} = \epsilon. 
    \end{align*}
    The optimal value of $p_{i}^{K}(x)$ on $X=[a,b]$ is $f_e^*$, see \eqref{eq:opt_selection}. It follows from Lemma~\ref{lem:value_distance} that $\left|f_e^*-f^*\right| \leq \epsilon$.

\subsection{Proof of Theorem~\ref{thm:alg_whole_complx}}\label{subsec:proof_complexity}
    In the stage of initialization, every agent $i$ evaluates $f_{i}(x)$ at $2m+1$ points to construct $p_{i}(x)$ of degree $m$. Hence, the orders of zeroth-order queries and flops are $\bigO{m}$ and $\bigO{mF_{0}}$, respectively. 
    The calculation in \eqref{eq:grid_point_evaluation} costs $\bigO{m}$ flops. For fixed $m_{i}$, through Fast Cosine Transform, we obtain the coefficients in \eqref{eq:coeff_value_relation} with $\bigO{m_{i} \log m_{i}}$ flops\cite{trefethen2013approximation}. Since evaluating $p_{i}(x)$ at one point requires $\bigO{m_{i}}$ flops via the Clenshaw algorithm\cite{boyd2014solving}, the check for \eqref{eq:adapt_stop_rule} costs $\bigOcom{m_{i}^{2}}$ flops. When $m_{i}$ is doubled from $2$ to $m$, the order of the total flops required by \eqref{eq:coeff_value_relation} and \eqref{eq:adapt_stop_rule} is of
    \begin{equation*}
        \bigO{m^{2}} + \bigOfrac{({\textstyle \frac{m}{2}})^{2}} + \ldots + \bigO{2^{2}} = \bigO{m^{2}}.
    \end{equation*}
    Therefore, the costs of this stage are $\bigO{m\cdot\max(m,F_{0})}$ flops and $\bigO{m}$ zeroth-order queries.

    Now we consider the second stage of iterations. Based on \eqref{eq:convergence_time}, we know that the order of the number of consensus iterations that guarantees \eqref{eq:itr_precision_requirement} is of
    \begin{equation*}
        U + \bigOfracEq{\log \frac{1}{\delta}} = \bigOfracEq{\log \frac{1}{\delta}} = \bigOfracEq{\log \frac{m}{\epsilon}}.
    \end{equation*}
    In every iteration, agent $i$ exchanges information with its neighbors only once and updates $p_{i}^{t}$ by \eqref{eq:avg_consensus_implement} with $\bigO{m}$ flops. Hence, the orders of communication and flops are of $\bigO{\log \frac{m}{\epsilon}}$ and $\bigO{m \log \frac{m}{\epsilon}}$, respectively.

    We discuss the final stage of polynomial optimization. The recurrence formula \eqref{eq:recurrence_coeff} requires $\bigO{m}$ flops. The eigenvalues of the colleague matrix $M_C$ are computed via the QR algorithm with $\mathcal{O}(m^{3})$ flops. The order of cost is brought down to $\mathcal{O}(m^{2})$ in Chebfun by recursion and exploitation of the special structure of $M_{C}$\cite{trefethen2013approximation}. When we search for the optimal value and optimal points in \eqref{eq:opt_selection}, the evaluation of $p_{i}^{K}(x)$ at one point requires $\bigO{m}$ flops by using the Clenshaw algorithm\cite{boyd2014solving}. Thus, \eqref{eq:opt_selection} costs $\mathcal{O}(m^{2})$ flops. Consequently, this stage needs $\mathcal{O}(m^{2})$ flops in total. We combine these results and arrive at Theorem~\ref{thm:alg_whole_complx}.

\subsection{Proof of Theorem~\ref{thm:sdp_equivalence}}\label{subsec:appendix_derivation}
The major fact that we use is
\begin{equation*}
    T_{u}(x)T_{v}(x) = \frac{1}{2} T_{u+v}(x) + \frac{1}{2} T_{|u-v|}(x),
\end{equation*}
which implies that the product of two Chebyshev polynomials equals a weighted sum of polynomials of certain degrees\cite{trefethen2013approximation}. Note that
\begin{equation}\label{eq:g_expansion}
    g_{i}^{K}(x) - t = (c_{0}'-t)T_{0}(x) + \sum_{j=1}^{m} c_{j}'T_{j}(x).
\end{equation}

When $m$ is odd, we have
\begin{equation*}
    \begin{split}
        &g_{i}^{K}(x) - t = (x+1)h_{1}^{2}(x) + (1-x)h_{2}^{2}(x) \\
             & = (T_{1}(x)+1)v_{1}(x)^{T}Qv_{1}(x) + (1-T_{1}(x))v_{2}(x)^{T}Q'v_{2}(x) \\
             & = \underbrace{\sum_{u,v} Q_{uv} T_{u}(x)T_{v}(x)}_{\numcircled{1}} + \underbrace{T_{1}(x) \sum_{u,v} Q_{uv} T_{u}(x)T_{v}(x)}_{\numcircled{2}} \\
             & \qquad + \underbrace{\sum_{u,v} Q'_{uv} T_{u}(x)T_{v}(x)}_{\numcircled{3}} - \underbrace{T_{1}(x)\sum_{u,v} Q'_{uv} T_{u}(x)T_{v}(x)}_{\numcircled{4}}.
    \end{split}
\end{equation*}
We first consider the coefficient of $T_{0}(x)$ in $g(x)$. From \eqref{eq:g_expansion}, this number is $c_{0}'-t$. For terms \numcircled{1} and \numcircled{3}, the added coefficient of $T_{0}(x)$ is $Q_{00} + Q'_{00} + \frac{1}{2} \big(\sum_{u} Q_{uu} + \sum_{u'} Q'_{u'u'} \big)$. For terms \numcircled{2} and \numcircled{4}, the added coefficient of $T_{0}(x)$ is $\sum_{|u-v|=1} \frac{1}{4} \left(Q_{uv} - Q'_{uv}\right)$. Hence,
\begin{equation}\label{eq:coeff_constraint_first_odd}
    \begin{split}
        c_{0}' - t &=  Q_{00} + Q'_{00} + \frac{1}{2} \Big(\sum_{u=1}^{d_{1}+1}Q_{uu} + \sum_{u=1}^{d_{2}+1} Q'_{uu}\Big) \\
                  &\quad + \frac{1}{4} \sum_{|u-v|=1} \left(Q_{uv} - Q'_{uv}\right).  
    \end{split}
\end{equation}
Then, we consider the coefficient of $T_{j}(x)$ in $g(x)$, where $j=1,\ldots,m$. From \eqref{eq:g_expansion}, it equals $c_{j}$. For terms \numcircled{1} and \numcircled{3}, when $(u,v) \in \mathcal{A}$ (see \eqref{eq:set_A}), 
the expansion of the product contains $T_{j}(x)$. Hence, the added coefficient of $T_{j}(x)$ is $\sum_{(u,v)\in \mathcal{A}} \frac{1}{2} \left(Q_{uv} + Q'_{uv}\right)$. For terms \numcircled{2} and \numcircled{4}, when $(u,v) \in \mathcal{B}$ (see \eqref{eq:set_B}), 
the expansion of the three-term product contains $T_{j}(x)$. Hence, the added coefficient of $T_{j}(x)$ is $\sum_{(u,v)\in \mathcal{B}} \frac{1}{4} \left(Q_{uv} - Q'_{uv}\right)$. It follows that for $j=1,\ldots,m$,
\begin{equation}\label{eq:coeff_constraint_general_odd}
        c_{j}' = \frac{1}{2} \sum_{(u,v)\in \mathcal{A}} \left(Q_{uv} + Q'_{uv} \right)  + \frac{1}{4} \sum_{(u,v)\in \mathcal{B}} \left(Q_{uv} - Q'_{uv}\right). 
\end{equation}
Equations~\eqref{eq:coeff_constraint_first_odd}, \eqref{eq:coeff_constraint_general_odd}, and the requirements that $Q \in \mathbb{S}_{+}^{d_{1} +1}, Q' \in \mathbb{S}_{+}^{d_{2} +1}$
constitute the constraints in \eqref{problem:odd_reformulated}. Therefore, when $m$ is odd, problem \eqref{problem:odd_reformulated} and problem \eqref{problem:poly_main} have the same optimal value.

When $m$ is even, we have
\begin{equation*}
   \begin{split}
       g_{i}^{K}&(x) - t = h_{1}^{2}(x) + (x+1)(1-x)h_{2}^{2}(x) \\
            &= v_{1}(x)^{T}Qv_{1}(x) + \frac{1-T_{2}(x)}{2} v_{2}(x)^{T}Q'v_{2}(x) \\
            &= \underbrace{\sum_{u,v} Q_{uv}T_{u}(x)T_{v}(x)}_{\numcircled{5}} + \underbrace{\frac{1}{2} \sum_{u,v} Q'_{uv}T_{u}(x)T_{v}(x)}_{\numcircled{6}} \\
            & \quad + \underbrace{\frac{1}{2} T_{2}(x) \sum_{u,v} Q'_{uv}T_{u}(x)T_{v}(x)}_{\numcircled{7}}.
   \end{split}
\end{equation*}
Similarly, we first consider the coefficient of $T_{0}(x)$ in $g(x)$, which equals $c_{0}'-t$. For terms \numcircled{5} and \numcircled{6}, the added coefficient of $T_{0}(x)$ is $Q_{00} + \frac{1}{2} Q'_{00} + \sum_{u} \frac{1}{2}Q_{uu} + \sum_{u'} \frac{1}{4}Q'_{u'u'}$. For term \numcircled{7}, the coefficient of $T_{0}(x)$ is $\frac{1}{2}\sum_{|u-v|=2} \frac{1}{4} Q'_{uv}$. Hence,
\begin{equation}\label{eq:coeff_constraint_first_even}
   \begin{split}
       c_{0}' - t &=  Q_{00} + \frac{1}{2} Q'_{00} + \frac{1}{2}\sum_{u=1}^{d_{1}+1}Q_{uu} \\
            &\quad + \frac{1}{4} \sum_{u=1}^{d_{2}+1} Q'_{uu} + \frac{1}{8} \sum_{|u-v|=2} Q'_{uv}.  
   \end{split}
\end{equation}
Then, we consider the coefficient of $T_{j}(x)$ in $g(x)$, where $j=1,\ldots,m$. From \eqref{eq:g_expansion}, it equals $c_{j}$. For terms \numcircled{5} and \numcircled{6}, when $(u,v) \in \mathcal{A}$ (see \eqref{eq:set_A}),
the expansion of the product contains $T_{j}(x)$. Hence, the added coefficient of $T_{j}(x)$ is $\sum_{(u,v)\in \mathcal{A}} \frac{1}{2} \left(Q_{uv} + Q'_{uv}\right)$. For term \numcircled{7}, when $(u,v) \in \mathcal{C} \triangleq \big\{(u,v)|u+v=i-2 \vee |u-v|=i-2 \vee |u+v-2|=i \vee \big||u-v|-2\big|=i \big\}$,
the expansion of the three-term product contains $T_{j}(x)$. Hence, the added coefficient of $T_{j}(x)$ is $\frac{1}{2} \sum_{(u,v)\in \mathcal{C}} \frac{1}{4} Q'_{uv}$. Consequently, for $j=1,\ldots,m$,
\begin{equation}\label{eq:coeff_constraint_general_even}
       c_{j}' = \frac{1}{2} \sum_{(u,v)\in \mathcal{A}} \Big(Q_{uv} + \frac{1}{2} Q'_{uv} \Big) + \frac{1}{8} \sum_{(u,v)\in \mathcal{C}} Q'_{uv}.
\end{equation}
Equations \eqref{eq:coeff_constraint_first_even}, \eqref{eq:coeff_constraint_general_even}, and the requirements of positive semi-definiteness $Q \in \mathbb{S}_{+}^{d_{1} +1}, Q' \in \mathbb{S}_{+}^{d_{2} +1}$ constitute the constraints in \eqref{problem:even_reformulated}. Therefore, when $m$ is even, problem \eqref{problem:even_reformulated} has the same optimal value as problem \eqref{problem:poly_main}.

    
    \balance
    \bibliographystyle{IEEEtran}
    \bibliography{journal}

\end{document}